\documentclass[12pt]{amsart}

\usepackage{amssymb,amsmath,array,latexsym,verbatim,color,mathtools}

\usepackage{subfigure}
\usepackage{tikz}
\usepackage{pgfmath}
\usepackage[hidelinks]{hyperref}
\usepackage{pgfplots}
\pgfplotsset{compat=1.15}
\usepackage{mathrsfs}
\usetikzlibrary{arrows}
\usetikzlibrary{shapes,arrows,calc,angles,positioning,intersections,quotes,decorations.markings,shapes.geometric,fit,datavisualization,datavisualization.formats.functions}
\usetikzlibrary{angles,positioning,intersections,quotes,decorations.markings,shapes.geometric,fit,datavisualization,datavisualization.formats.functions}
\usetikzlibrary{matrix, arrows.meta}
\usetikzlibrary{automata,shapes}
\tikzset{
  invisible/.style={opacity=0},
  visible on/.style={alt=#1{}{invisible}},
  alt/.code args={<#1>#2#3}{%
    \alt<#1>{\pgfkeysalso{#2}}{\pgfkeysalso{#3}} 
  },
  state st/.style={draw,circle,top color=orange!2,bottom color=red!50!orange!50}
}

\providecommand{\pgfsyspdfmark}[3]{}
\newcommand{\la}{\langle \, }
\newcommand{\ra}{\, \rangle}
\newcommand{\sm}{\mathcal S}

\newcommand{\srg}{\mathrm{SRG}}
\newcommand{\tr}{\mathrm{tr}}
\newcommand{\msr}{\mathrm{mr}_+}

\newcommand{\C}{\operatorname{\mathsf{C}}} 

\newtheorem{thm}{Theorem}[section] 
\newtheorem{lem}[thm]{Lemma}   
\newtheorem{cor}[thm]{Corollary} 
\newtheorem{prop}[thm]{Proposition}

\theoremstyle{definition}
\newtheorem{defn}[thm]{Definition}

\title[On distinct eigenvalues of triangle-free SRGs]{On the minimum number of distinct eigenvalues of triangle-free strongly regular graphs}
\author[E. Egolf]{Emily~Egolf}
\address{Emily Egolf, Department of Mathematics, Fort Lewis College, Durango, CO  81301, USA}
\email{eaegolf@fortlewis.edu}
\author[V. Furst]{Veronika~Furst}
\address{Veronika Furst, Department of Mathematics, Fort Lewis College, Durango, CO  81301, USA}
\email{furst\_v@fortlewis.edu}

\begin{document}

\begin{abstract}
    Among the seven known (non-degenerate) triangle-free strongly regular graphs, we prove that the Clebsch graph describes a matrix with exactly two distinct eigenvalues while five of the graphs do not.  In showing that the minimum number of distinct eigenvalues of the Sims-Gewirtz graph is three, we answer a recently stated open question.
\end{abstract}

\maketitle



\section{Introduction}

The {\em inverse eigenvalue problem of a graph (IEPG)}, in its most broad sense, is to determine all possible spectra of real symmetric matrices described by a graph.  That is, we consider all matrices with off-diagonal zero/non-zero pattern matching that of the graph's adjacency matrix.  Several subproblems of the IEPG have gained much attention over the past decade, in particular, the minimum number of distinct eigenvalues of a graph, and more specifically, when a graph has a corresponding matrix that has exactly two distinct eigenvalues.  (The spectral decomposition shows that a graph that admits a matrix with one eigenvalue must be edgeless.)  For background on the IEPG, see the book \cite{MR4478249} and the references therein; results on graphs that admit two distinct eigenvalues can be found in \cite{MR3231825, MR3118943, MR4044603, MR3781593, MR4214542, MR4603836, MR4646344, MR3891770, MR4466746}.  

For standard graph theory terminology and notation, see e.g. \cite{MR1367739}.  Given a (simple, undirected) graph $G$, we let $V(G)$ and $E(G)$ represent its set of vertices and edges, respectively.  If vertices $u$ and $v$ are adjacent, we write $u\sim v$ (equivalently, $v\sim u$), and let $N(u) = \{v\in V(G): \ v\sim u\}$ be the open neighborhood of $u$.  We denote the edge $\{u,v\}$ as $(u,v)$ and the cycle $C_k$ on $k$ vertices $v_1, \ldots, v_k$ as $(v_1, \ldots, v_k, v_1)$.  In Section \ref{SimsGewirtz} we will also write $C_k$ in terms of its edges as $[(v_1,v_2), (v_2,v_3), \ldots, (v_{k-1},v_k)]$.  The complete graph and complete multipartite graph are $K_n$ and $K_{n_1,n_2,\ldots,n_k}$, respectively.  For $S\subseteq V(G)$, let $G[S]$ be the subgraph of $G$ induced by $S$; let $G-S = G[V(G)\setminus S]$, and if $H$ is a subgraph of $G$, let $G-H = G - V(H)$.  Since we will be interested in vector representations of graphs, write $\langle f,g \rangle$ for the standard dot product of two real vectors $f,g \in \mathbb{R}^d$.

To a graph $G$ on $n$ vertices, we associate the set of $n\times n$ real symmetric matrices described by $G$:
\[ \sm(G) =  \{ A = [a_{ij}] = A^T: \ a_{ij} = 0 \mbox{ for } i\neq j \mbox{ if and only if } (v_i,v_j) \notin E(G) \}.
\]
No restrictions are placed on the diagonal entries of the matrix $A$ while the off-diagonal zero/non-zero entries of $A$ appear in the same locations as those of the adjacency matrix of the graph.  Let $q(A)$ denote the number of distinct eigenvalues of a square matrix $A$, and define the {\em minimum number of distinct eigenvalues} of $G$ to be 
\[ q(G) = \min \{q(A): \ A \in \sm(G) \}. \]

A graph $G$ is said to be {\em strongly regular} with (nonnegative) parameters
$(n, k, \lambda, \mu)$ if $|V(G)| = n$, each vertex of $G$ has degree $k$
(so $G$ is $k$-regular), any two adjacent vertices have exactly $\lambda$ neighbors
in common, and any two non-adjacent vertices have exactly $\mu$ neighbors in common.
We abbreviate the term strongly regular graph as SRG. It is customary to
omit the degenerate cases when $k = n-1$, $\lambda = k-1$, or $\mu = 0$,
which correspond to disjoint unions of $(k+1)$-cliques, and their complements, the complete multipartite graphs with partite sets of order $k+1$.  The adjacency matrix of any (non-degenerate and hence connected) strongly regular graph has three distinct eigenvalues \cite{MR1367739}, so if $G$ is a strongly regular graph, then $q(G) \in \{2,3\}$.  (It is well known that
$q(K_n) = q(K_{m,m,\ldots,m}) = 2$ \cite{MR3118943}.)  

When $\lambda = 0$, the SRG is known as {\em triangle-free}.
Not including the complete bipartite graphs $K_{m,m}$, there are seven known
triangle-free strongly regular graphs (see \cite{aebrouwer}):
\begin{itemize}
	\item Pentagon $C_5 = \srg(5,2,0,1)$
	\item Petersen graph $\srg(10,3,0,1)$
	\item Clebsch graph $\srg(16,5,0,2)$
	\item Hoffman-Singleton graph $\srg(50,7,0,1)$
	\item Sims-Gewirtz graph $\srg(56,10,0,2)$
	\item Mesner-M22 graph $\srg(77,16,0,4)$
	\item Higman-Sims graph $\srg(100,22,0,6)$
\end{itemize}

In this paper, when we refer to triangle-free SRGs, we assume that these seven comprise the complete list.  We find the $q$-values for six of the seven.  Section \ref{prelim} contains necessary background results, from which $q=3$ can easily be established for four of the seven triangle-free SRGs.  In Section \ref{Clebsch} we prove $q=2$ for the Clebsch graph; although this result appears in the recent preprint \cite{fallat2024}, our technique is different and introduces notation and a key result that will be used in Section \ref{SimsGewirtz}.  Section 4 contains our main result, namely, that $q=3$ for the Sims-Gewirtz graph; this answers an open question in \cite{fallat2024}.  We leave the same question open for the Higman-Sims graph.

\section{Preliminary Results} \label{prelim}

For a graph $G$, we consider the subclass of positive semidefinite matrices in $\sm(G)$.  The {\em minimum positive semidefinite rank} of a graph $G$ is defined to be $\msr(G) = \min \{\mathrm{rank}(A): \ A\in \sm(G) \mbox{ and $A$ is positive semidefinite} \}$ (see \cite{MR4478249}). 
Every graph $G$ has a {\em faithful orthogonal representation} in $\mathbb{R}^d$ where $d\geq \msr(G)$ \cite{MR514926}. 
That is, let
$\{f_1, \ldots, f_n\}\subseteq \mathbb{R}^d$ be a spanning set such that there is a one-to-one
correspondence between the vertices $v_i \in V(G)$ and the vectors $f_i$ with $v_i \sim v_j$
if and only if $\langle f_i, f_j \rangle \neq 0$. Consider the {\em Gram matrix}
$M = F^TF$ where $F = [ f_1, \ldots, f_n ]$ is the $d\times n$ matrix whose
columns are the vectors in the orthogonal representation of $G$. Then $M_{ij} = \langle f_i, f_j \rangle$, and it follows that $M \in \sm(G)$.

A graph $G$ with at least one edge satisfies $q(G) = 2$ if and only if
$\sm(G)$ contains a symmetric orthogonal matrix \cite{MR3118943}.  By Lemma 2.3 of \cite{MR3118943}, we can shift the eigenvalues of this matrix from $\{-1,1\}$ to $\{0,1\}$ to yield a positive semidefinite matrix, which can be
decomposed, using the spectral theorem, as a Gram matrix $M = F^TF$ where
$FF^T = I_d$. Note that $M^2 = M$.  In fact, $q(G) = 2$ if and only if $\sm(G)$ contains a Gram matrix $M = F^TF = [ \langle f_i, f_j \rangle ]$ such that $M^2 = M$ (see \cite{MR3781593}, \cite{MR4214542}).

The following result (a corollary of Theorem 3.2 of \cite{MR3118943}) is a
very convenient necessary, though not sufficient, condition for a graph $G$ to
satisfy $q(G)=2$. For completeness we include the proof.

\begin{lem}
	If $q(G) = 2$, then any pair of non-adjacent vertices in $G$ cannot have a
	unique common neighbor.
\end{lem}

\begin{proof}
	Let $M \in \sm(G)$ such that $M = [ \langle f_i, f_j \rangle ]$ and $M^2 = M$.
	If $v_k$ is the unique common neighbor of non-adjacent vertices $v_i$ and $v_j$,
	then 
	\begin{align*}
		\langle f_i, f_j \rangle = M_{ij} = M^2_{ij} =
		\sum_{l} \langle f_i, f_l \rangle \langle f_l, f_j \rangle =
		\langle f_i, f_k \rangle \langle f_k, f_j \rangle.
	\end{align*}
	The contradiction follows by noting that the left-hand side is zero while the
	right-hand side is nonzero.
\end{proof}

\begin{cor}
	The Pentagon, Petersen graph, and Hoffman-Singleton graph have $q=3$.
\end{cor}

\begin{prop} \label{Gramdiag}
	Let $G = \srg(n, k, 0, \mu)$ and suppose $q(G)=2$. Let $M \in
	\sm(G)$ be such that $M = [ \langle f_i, f_j \rangle ]$ and $M^2 = M$.  Then every
	diagonal entry of $M$ must be $\frac{1}{2}$.
\end{prop}

\begin{proof}
	If $v_i \sim v_j$, then
	\begin{align*}
		\langle f_i, f_j \rangle = M_{ij} = M^2_{ij} =
		\sum_l \langle f_i, f_l \rangle \langle f_l, f_j \rangle =
		\langle f_i, f_i \rangle \langle f_i, f_j \rangle + \langle f_i, f_j \rangle \langle f_j, f_j \rangle
	\end{align*}
	since $v_i$ and $v_j$ share no common neighbors. Since
	$\langle f_i, f_j \rangle \neq 0$, we get
	\begin{align*}
		1 = \langle f_i, f_i \rangle + \langle f_j, f_j \rangle
		= \|f_i\|^2 + \|f_j\|^2.
	\end{align*}
	If $v_{i_0}$ belongs to an odd cycle, then the vectors corresponding to the
	vertices of the cycle must alternate in length between $\|f_{i_0}\|^2$ and
	$1 - \|f_{i_0}\|^2$, which implies $\|f_{i_0}\|^2 = \frac{1}{2}$. Since every
	vertex of $G$ belongs to a 5-cycle (verified by noting that the diagonal
	entries of the fifth power of each adjacency matrix are nonzero, and, in fact, equal, while $G$ contains no 3-cycles), we see that
	$M_{ii} = \|f_i\|^2 = \frac{1}{2}$ for every $i$.
\end{proof}

\begin{cor}
	If $G = \srg(n, k, 0, \mu)$ with $q(G)=2$, then $n$ is even.
\end{cor}

\begin{proof}
	By Proposition \ref{Gramdiag}, we know there exists Gram matrix $M \in \sm(G)$
	such that $M = F^TF$, $I_d = FF^T$, and
	$\mathrm{diag}(M) = (\frac{1}{2}, \frac{1}{2}, \ldots, \frac{1}{2})$. By the commutativity of the trace of a matrix product,
	\[ d = \tr(FF^T) = \tr(F^TF) = \frac{n}{2},\]
	so $n = 2d$.
\end{proof}

\begin{cor}
	The Mesner-M22 graph has $q=3$.
\end{cor}

The remainder of this paper is dedicated to establishing the $q$-values of the
Clebsch and Sims-Gewirtz graphs.

\section{The Clebsch Graph} \label{Clebsch}

\subsection{The Clebsch Graph Describes a Matrix with Two Distinct Eigenvalues} \label{q(Clebsch)=2}

The following simple lemma is crucial in establishing both that $q=2$ for the Clebsch graph and $q\neq 2$ for the Sims-Gewirtz graph.

\begin{lem} \label{signlemma}
	Let $G = \srg(n,k,0,2)$ with $q(G) = 2$, and let $M \in \sm(G)$ such that
	$M = [ \la f_i, f_j \ra ]$ and $M^2 = M$. If $(v_i, v_j, v_k, v_l, v_i)$ is a
	4-cycle in $G$, then
	\begin{align*}
		\langle f_j, f_k \rangle = \pm \langle f_i, f_l \rangle\quad\mbox{and}\quad \langle f_i, f_j \rangle = \mp \langle f_k, f_l \rangle
	\end{align*}
	where exactly one choice is positive and one choice is negative.
\end{lem}

\begin{proof}
	Since $v_i$ and $v_k$ are non-adjacent vertices, and $v_j$ and $v_l$ are their
	only common neighbors, we have
	\begin{equation}\label{firsteq}
		0 = M_{ik} = M^2_{ik} = \langle f_i, f_j \rangle \langle f_j, f_k \rangle + \langle f_i, f_l \rangle \langle f_l, f_k \rangle
	\end{equation}
	where every term on the right-hand side is nonzero.  Similarly,
	\begin{equation}\label{secondeq}
		0 = M_{jl} = M^2_{jl} = \langle f_j, f_k \rangle \langle f_k, f_l \rangle + \langle f_j, f_i \rangle \langle f_i, f_l \rangle.
	\end{equation}
	Solving Equation \ref{secondeq} for $\langle f_i, f_j \rangle$ and plugging
	into Equation \ref{firsteq}, we obtain
	\begin{align*}
		0 = -\dfrac{\langle f_j, f_k \rangle \langle f_k, f_l \rangle}{\langle f_i, f_l \rangle} \langle f_j, f_k \rangle + \langle f_i, f_l \rangle \langle f_l, f_k \rangle,
	\end{align*}
	which implies
	\begin{align*}
		\langle f_j, f_k \rangle^2 = \langle f_i, f_l \rangle^2
	\end{align*}
	Similarly, solving Equation \ref{secondeq} for $\langle f_j, f_k \rangle$ and
	plugging into Equation \ref{firsteq} yields
	\begin{align*}
		\langle f_i, f_j \rangle^2 = \langle f_k, f_l \rangle^2
	\end{align*}
	The conclusion follows by Equation \ref{firsteq}.
\end{proof}

\begin{figure}[h!]
	\begin{tikzpicture}[scale=.7,line cap=round,line join=round,>=triangle 45,x=1cm,y=1cm,V/.style={circle, inner sep = 2pt,fill=black}]
\draw [line width=2pt] (-3.8505784485542613,5.028958908075179)-- (-3.8505784485542613,-3.9710410919248122);
\draw [line width=2pt] (-3.8505784485542613,-3.9710410919248122)-- (6.149421551445739,-3.9710410919248122);
\draw [line width=2pt] (6.149421551445739,-3.9710410919248122)-- (6.149421551445739,5.028958908075179);
\draw [line width=2pt] (6.149421551445739,5.028958908075179)-- (-3.8505784485542613,5.028958908075179);
\draw [line width=2pt] (6.149421551445739,5.028958908075179)-- (4.149421551445739,4.028958908075179);
\draw [line width=2pt] (4.149421551445739,4.028958908075179)-- (4.149421551445739,-4.971041091924821);
\draw [line width=2pt] (4.149421551445739,-4.971041091924821)-- (6.149421551445739,-3.9710410919248122);
\draw [line width=2pt] (-3.8505784485542613,-3.9710410919248122)-- (-5.850578448554261,-4.971041091924821);
\draw [line width=2pt] (-5.850578448554261,-4.971041091924821)-- (4.149421551445739,-4.971041091924821);
\draw [line width=2pt] (4.149421551445739,4.028958908075179)-- (-5.850578448554261,4.028958908075179);
\draw [line width=2pt] (-5.850578448554261,4.028958908075179)-- (-5.850578448554261,-4.971041091924821);
\draw [line width=2pt] (-5.850578448554261,4.028958908075179)-- (-3.8505784485542613,5.028958908075179);
\draw [line width=2pt] (-1,2)-- (-1,-1.5);
\draw [line width=2pt] (-1,-1.5)-- (2.5,-1.5);
\draw [line width=2pt] (2.5,-1.5)-- (2.5,2);
\draw [line width=2pt] (2.5,2)-- (-1,2);
\draw [line width=2pt] (2.5,2)-- (1.5,1.5);
\draw [line width=2pt] (1.5,1.5)-- (1.5,-2);
\draw [line width=2pt] (1.5,-2)-- (2.5,-1.5);
\draw [line width=2pt] (-1,-1.5)-- (-2,-2);
\draw [line width=2pt] (-2,-2)-- (1.5,-2);
\draw [line width=2pt] (1.5,1.5)-- (-2,1.5);
\draw [line width=2pt] (-2,1.5)-- (-2,-2);
\draw [line width=2pt] (-2,1.5)-- (-1,2);
\draw [line width=2pt] (-3.8505784485542613,5.028958908075179)-- (-1,2);
\draw [line width=2pt] (-2,1.5)-- (-5.850578448554261,4.028958908075179);
\draw [line width=2pt] (-3.8505784485542613,-3.9710410919248122)-- (-1,-1.5);
\draw [line width=2pt] (-2,-2)-- (-5.850578448554261,-4.971041091924821);
\draw [line width=2pt] (1.5,-2)-- (4.149421551445739,-4.971041091924821);
\draw [line width=2pt] (2.5,-1.5)-- (6.149421551445739,-3.9710410919248122);
\draw [line width=2pt] (2.5,2)-- (6.149421551445739,5.028958908075179);
\draw [line width=2pt] (4.149421551445739,4.028958908075179)-- (1.5,1.5);

\begin{scope}[nodes=V]
	\node [fill=black,label={below right:0100}] at (-3.8505784485542613,-3.9710410919248122) {};
	\node [fill=black,label={below right:0101}] at (6.149421551445739,-3.9710410919248122) {};
	\node [fill=black,label={above:1101}] at (6.149421551445739,5.028958908075179) {};
	\node [fill=black,label={above:1100}] at (-3.8505784485542613,5.028958908075179) {};
	\node [fill=black,label={above left:1000}] at (-5.850578448554261,4.028958908075179) {};
	\node [fill=black,label={below left:0000}] at (-5.850578448554261,-4.971041091924821) {};
	\node [fill=black,label={below right:0001}] at (4.149421551445739,-4.971041091924821) {};
	\node [fill=black,label={above left:1001}] at (4.149421551445739,4.028958908075179) {};
	\node [fill=black,label={above right:0110}] at (-1,-1.5) {};
	\node [fill=black,label={above right:0111}] at (2.5,-1.5) {};
	\node [fill=black,label={below right:1111}] at (2.5,2) {};
	\node [fill=black,label={above right:1110}] at (-1,2) {};
	\node [fill=black,label={below left:1010}] at (-2,1.5) {};
	\node [fill=black,label={above left:0010}] at (-2,-2) {};
	\node [fill=black,label={below left:0011}] at (1.5,-2) {};
	\node [fill=black,label={below left:1011}] at (1.5,1.5) {};
\end{scope}
\end{tikzpicture}
	\caption{Clebsch graph shown as the 4-cube missing eight additional edges $(v,-v)$ for $v = 0000, 0001, 0010, 0011, 0100, 0101, 0110, 0111$.}
    \label{fig:Clebsch}
\end{figure}
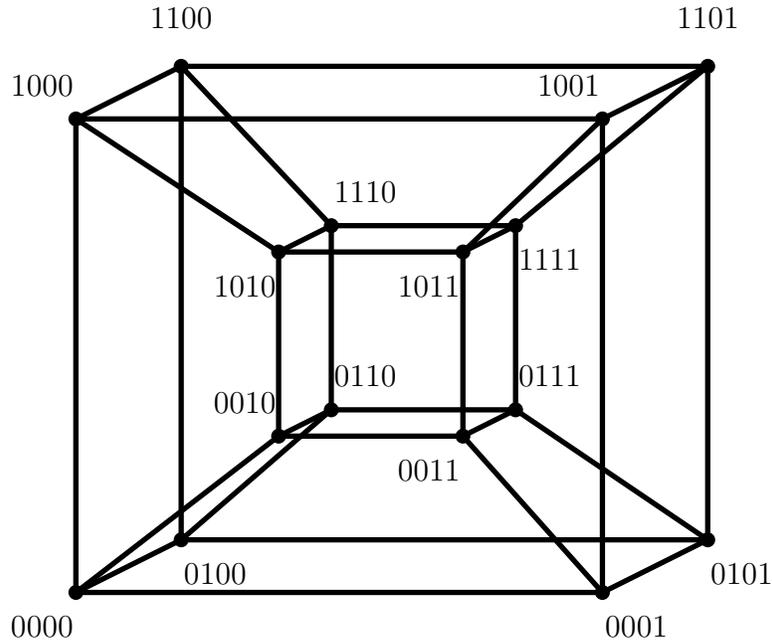

In what follows, we denote the Clebsch graph $\srg(16,5,0,2)$ by $\C$.  One way to describe $\C$ is to label the vertices of the 4-cube with 4-bit binary numbers, and two vertices are adjacent if and only if their labels differ in exactly one bit or in all four bits.  Figure \ref{fig:Clebsch} shows $\C$ without the eight edges that connect a vertex $v$ to its binary complement.

To construct a matrix $M\in \sm(\C)$ with two distinct eigenvalues, 
consider the five edges incident to vertex 0000 and denote the corresponding
entries of $M$ by $a,b,c,d,e$. By Lemma \ref{signlemma}, every nonzero
off-diagonal entry of $M$ is $\pm a, \pm b, \pm c, \pm d, \pm e$. In fact, it
turns out that we may take $a = b = c = d = e = \frac{1}{2\sqrt{5}}$ to obtain
the appropriate matrix.

\begin{thm} \label{thm:Clebsch-q2}
    For the Clebsch graph $\C$, $q(\C) = 2$.
\end{thm}

\begin{proof}
    The matrix
\begin{small}
	\[ M = \frac{1}{2\sqrt{5}} \begin{bmatrix}
		\sqrt{5} & -1 & 1 & 0 & 1 & 0 & 0  & 0 & 1 & 0 & 0 & 0 & 0 & 0 & 0 & 1 \\
		-1 & \sqrt{5} & 0 & 1 & 0 & 1 & 0 & 0 & 0 & 1 & 0 & 0 & 0 & 0 & -1 & 0 \\
		1 & 0 & \sqrt{5} & 1 & 0 & 0 & 1 & 0 & 0 & 0 & -1 & 0 & 0 & -1 & 0 & 0 \\
		0 & 1 & 1 & \sqrt{5} & 0 & 0 & 0 & -1 & 0 & 0 & 0 & 1 & 1 & 0 & 0 & 0 \\
		1 & 0 & 0 & 0 & \sqrt{5} & 1 & -1 & 0 & 0 & 0 & 0 & -1 & 1 & 0 & 0 & 0 \\
		0 & 1 & 0 & 0 & 1 & \sqrt{5} & 0 & 1 & 0 & 0 & 1 & 0 & 0 & -1 & 0 & 0 \\
		0 & 0 & 1 & 0 & -1 & 0 & \sqrt{5} & 1 & 0 & 1 & 0 & 0 & 0 & 0 & 1 & 0 \\
		0 & 0 & 0 & -1 & 0 & 1 & 1 & \sqrt{5} & -1 & 0 & 0 & 0 & 0 & 0 & 0 & 1 \\
		1 & 0 & 0 & 0 & 0 & 0 & 0 & -1 & \sqrt{5} & 1 & 1 & 0 & -1 & 0 & 0  & 0 \\
		0 & 1 & 0 & 0 & 0 & 0 & 1 & 0 & 1 & \sqrt{5} & 0 & -1 & 0 & 1 & 0 & 0 \\
		0 & 0 & -1 & 0 & 0 & 1 & 0 & 0 & 1 & 0 & \sqrt{5} & 1 & 0 & 0 & 1 & 0 \\
		0 & 0 & 0 & 1 & -1 & 0 & 0 & 0 & 0 & -1 & 1 & \sqrt{5} & 0 & 0 & 0 & 1 \\
		0 & 0 & 0 & 1 & 1 & 0 & 0 & 0 & -1 & 0 & 0 & 0 & \sqrt{5} & 1 & 1 & 0 \\
		0 & 0 & -1 & 0 & 0 & -1 & 0 & 0 & 0 & 1 & 0 & 0 & 1 & \sqrt{5} & 0 & 1 \\
		0 & -1 & 0 & 0 & 0 & 0 & 1 & 0 & 0 & 0 & 1 & 0 & 1 & 0 & \sqrt{5} & -1 \\
		1 & 0 & 0 & 0 & 0 & 0 & 0 & 1 & 0 & 0 & 0 & 1 & 0 & 1 & -1 & \sqrt{5} \\
	\end{bmatrix} \]
\end{small}
has exactly two
eigenvalues $0$ (with multiplicity 8) and $1$ (with multiplicity 8), and it is easy to verify that $M \in \sm(\C)$.
\end{proof}

\begin{cor}
    For the Clebsch graph $\C$, $\msr(\C) = 8$.
\end{cor}

\begin{proof}
For a graph $G$ on $n$ vertices with $q(G) = 2$, it is known that $\msr(G) \leq \lfloor \frac{n}{2} \rfloor$ \cite{MR4044603, MR4214542}, so $\msr(\C) \leq 8$.  On the
other hand, a triangle-free graph $G$ on $n$ vertices with no isolated
vertices satisfies $\msr(G) \geq \frac{n}{2}$ \cite{MR1116372} (see also
\cite{MR2775783}). So $\msr(\C) = \frac{16}{2} = 8$ and the Clebsch graph is a Rosenfeld graph.  The positive semidefinite matrix $M$ constructed in the previous proof can be decomposed as $M = F^TF$ where $FF^T$ is the $8\times 8$ identity matrix.
\end{proof}

\subsection{The Plus Graph of a Graph} The emphasis on 4-cycles, due to Lemma \ref{signlemma}, motivates the next definition.

\begin{defn} \label{defn:plusgraph}
Given a graph $G$, define its {\em plus graph} $G^+$ as the graph with vertex
set $V(G^+) = \{e^+: \ e\in E(G)\}$ and edge set
	\begin{align*}
	E(G^+) &= \{e^+f^+: \ \mbox{$e$ and $f$ are a pair of non-incident edges of some $4$-cycle in $G$} \} \\
    &= \{e^+f^+: \ \mbox{$\left\{e,f\right\}$ is a perfect matching of some $4$-cycle in $G$} \}.
	\end{align*}
\end{defn}

\begin{figure}[h!]
    \centering
    \begin{tikzpicture}[line cap=round,line join=round,>=triangle 45,x=1cm,y=1cm,V/.style={circle, inner sep = 2pt,fill=black}]
\draw [line width=2pt] (1,0)-- (0,1) node [midway] {};
\draw [line width=2pt] (0,1)-- (-1,0) node [midway] {};
\draw [line width=2pt] (-1,0)-- (0,-1) node [midway] {};
\draw [line width=2pt] (0,-1)-- (1,0) node [midway] {};
\draw [line width=2pt] (-1,0)-- (1,0) node [midway] {};

\begin{scope}[nodes=V]
\node at (1,0) {};
\node at (0,1) {};
\node at (-1,0) {};
\node at (0,-1) {};
\end{scope}

\draw [line width=2pt] (7,0)-- (6,1) node [midway] {};
\draw [line width=2pt] (6,1)-- (5,0) node [midway] {};
\draw [line width=2pt] (5,0)-- (6,-1) node [midway] {};
\draw [line width=2pt] (6,-1)-- (7,0) node [midway] {};
\draw [line width=2pt] (8.44,0)-- (7,0) node [midway] {};

\begin{scope}[nodes=V]
\node at (7,0) {};
\node at (6,1) {};
\node at (5,0) {};
\node at (6,-1) {};
\node at (8.44,0) {};
\end{scope}
\end{tikzpicture}
    \caption{The diamond graph and the banner graph.}
    \label{fig:plusgraphs}
\end{figure}
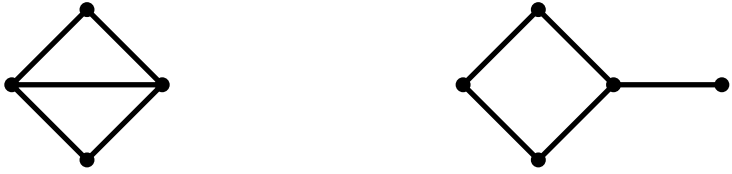

Note that $C_4^+ \cong K_2 \cup K_2$, which can be visualized as a plus sign.  The map $G \mapsto G^+$ is well-defined since any pair of edges $(x,y)$ and $(z,w)$ are non-incident edges of a $4$-cycle in $G$ if and only $(\varphi(x),\varphi(y))$ and $(\varphi(z),\varphi(w))$ are non-incident edges of a $4$-cycle for any graph isomorphism $\varphi$.  However, this map is clearly not one-to-one:  for example, $K_2\cup K_2 \cup K_1$ is the plus graph of both the diamond graph and the banner graph (see Figure \ref{fig:plusgraphs}).

In fact, we can relax the previous definition by starting with a subset $S$ of edges in a graph $G$ and defining the graph $S^{+_G}$ by $V(S^{+_G}) = \{e^+: \ e\in S\}$ and edge set $E(S^{+_G})$ as in Definition \ref{defn:plusgraph}.  For a subgraph $H$ of $G$, we will write $H^{+_G}$ for $E(H)^{+_G}$.  Note that $G^+ = E(G)^{+_G}$.

\begin{lem} \label{lem:plus-subgraph}
    If $G$ is a graph and $S \subseteq E(G)$, then $S^{+_G}$ is a subgraph of $G^+$.
\end{lem}

The following construction appears in the form of the 3-cube in the Clebsch graph and will become critical in Section \ref{SimsGewirtz} for the Sims-Gewirtz graph.

\begin{defn} \label{defn:trapgph}
An $n$-{\em trapezohedral graph} $T_n$ is a graph on $2n+2$ vertices in which $2n$ vertices comprise an outer cycle $C_{2n} = (a_0, b_0, a_1, b_1, \ldots, a_{n-1}, b_{n-1},a_0)$ and the two remaining vertices $\alpha, \beta$ satisfy $N(\alpha) = \{a_0, \ldots, a_{n-1}\}$ and $N(\beta) = \{b_0, \ldots, b_{n-1}\}$.
\end{defn}

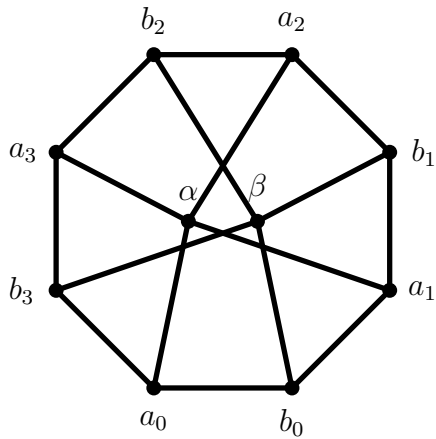
\begin{figure}[h]
	\definecolor{uuuuuu}{rgb}{0.26666666666666666,0.26666666666666666,0.26666666666666666}
\begin{tikzpicture}[scale=.4,line cap=round,line join=round,>=triangle 45,
	x=1cm,y=1cm, every node/.style={circle,inner sep = 2pt,fill=black}]
\draw [line width=2pt] (-2.2961005941905386,5.54327719506772)-- (2.2961005941905386,5.54327719506772);
\draw [line width=2pt] (2.2961005941905386,5.54327719506772)-- (5.54327719506772,2.2961005941905386);
\draw [line width=2pt] (5.54327719506772,2.2961005941905386)-- (5.54327719506772,-2.2961005941905386);
\draw [line width=2pt] (5.54327719506772,-2.2961005941905386)-- (2.2961005941905386,-5.54327719506772);
\draw [line width=2pt] (2.2961005941905386,-5.54327719506772)-- (-2.2961005941905386,-5.54327719506772);
\draw [line width=2pt] (-2.2961005941905386,-5.54327719506772)-- (-5.54327719506772,-2.2961005941905386);
\draw [line width=2pt] (-5.54327719506772,-2.2961005941905386)-- (-5.54327719506772,2.2961005941905386);
\draw [line width=2pt] (-5.54327719506772,2.2961005941905386)-- (-2.2961005941905386,5.54327719506772);
\draw [line width=2pt] (-1.1480502970952693,0)-- (-2.2961005941905386,-5.54327719506772);
\draw [line width=2pt] (-1.1480502970952693,0)-- (-5.54327719506772,2.2961005941905386);
\draw [line width=2pt] (-1.1480502970952693,0)-- (2.2961005941905386,5.54327719506772);
\draw [line width=2pt] (-1.1480502970952693,0)-- (5.54327719506772,-2.2961005941905386);
\draw [line width=2pt] (1.1480502970952693,0)-- (-5.54327719506772,-2.2961005941905386);
\draw [line width=2pt] (1.1480502970952693,0)-- (2.2961005941905386,-5.54327719506772);
\draw [line width=2pt] (1.1480502970952693,0)-- (5.54327719506772,2.2961005941905386);
\draw [line width=2pt] (1.1480502970952693,0)-- (-2.2961005941905386,5.54327719506772);
	\node [label={below:$b_0$}] at (2.2961005941905386,-5.54327719506772) {};
	\node [label={above:$b_2$}] at (-2.2961005941905386,5.54327719506772) {};
	\node [label={below:$a_0$}] at (-2.2961005941905386,-5.54327719506772) {};
	\node [label={above:$a_2$}] at (2.2961005941905386,5.54327719506772) {};
	\node [label={left:$b_3$}] at (-5.54327719506772,-2.2961005941905386) {};
	\node [label={right:$b_1$}] at (5.54327719506772,2.2961005941905386) {};
	\node [label={right:$a_1$}] at (5.54327719506772,-2.2961005941905386) {};
	\node [label={left:$a_3$}] at (-5.54327719506772,2.2961005941905386) {};
	\node [label={above:$\alpha$}] at (-1.1480502970952693,0) {};
	\node [label={above:$\beta$}] at (1.1480502970952693,0) {};
\end{tikzpicture}
    \caption{The trapezohedral graph $T_4$.}
    \label{fig:tsg}
\end{figure}

\begin{thm} \label{tsgConn}
	Let $T_n$ be the trapezohedral graph on $2n+2$ vertices.  Then $T_n^+$ is connected if and only if $n\pmod 3 \not\equiv 0$.
\end{thm}

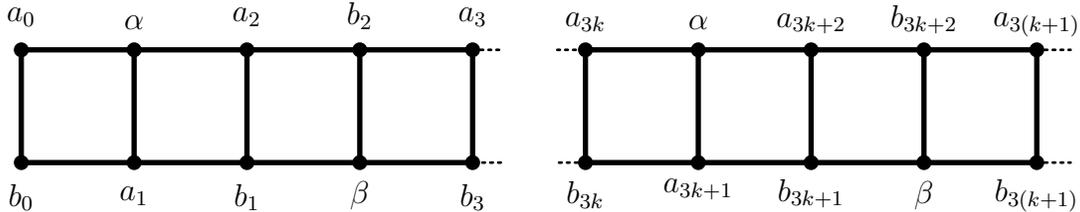
\begin{figure}[h!]
    \centering
    \begin{tikzpicture}[scale=.75,line cap=round,line join=round,>=triangle
	45,x=1cm,y=1cm,V/.style={circle, inner sep = 2pt,fill=black}]
	\draw [line width=2pt] (-8,2)-- (-8,0);
	\draw [line width=2pt] (8,2)-- (8,0);
	\draw [line width=2pt] (6,2)-- (6,0);
	\draw [line width=2pt] (4,2)-- (4,0);
	\draw [line width=2pt] (2,2)-- (2,0);
	\draw [line width=2pt] (0,2)-- (0,0);
	\draw [line width=2pt] (-2,2)-- (-2,0);
	\draw [line width=2pt] (-4,2)-- (-4,0);
	\draw [line width=2pt] (-6,2)-- (-6,0);
	\draw [line width=2pt] (-8,2)-- (0,2);
	\draw [line width=2pt] (-8,0)-- (0,0);
	\draw [line width=1pt,dotted] (0,2)-- (0.5,2);
	\draw [line width=1pt,dotted] (1.5,2)-- (2,2);
	\draw [line width=1pt,dotted] (2,0)-- (1.5,0);
	\draw [line width=1pt,dotted] (0.5,0)-- (0,0);
	\draw [line width=2pt] (2,2)-- (8,2);
	\draw [line width=2pt] (2,0)-- (8,0);
	\draw [line width=2pt] (8,2)-- (10,2);
	\draw [line width=2pt] (10,2)-- (10,0);
	\draw [line width=2pt] (10,0)-- (8,0);
	\draw [line width=1pt,dotted] (10,2)-- (10.6,2);
	\draw [line width=1pt,dotted] (10,0)-- (10.6,0);
	\begin{scope}[nodes=V]
		\node [label={below:$b_0$}] at (-8,0) {};
		\node [label={above:$a_0$}] at (-8,2) {};
		\node [label={below:$a_1$}] at (-6,0) {};
		\node [label={above:$\alpha$}] at (-6,2) {};
		\node [label={below:$b_1$}] at (-4,0) {};
		\node [label={above:$a_2$}] at (-4,2) {};
		\node [label={below:$\beta$}] at (-2,0) {};
		\node [label={above:$b_2$}] at (-2,2) {};
		\node [label={below:$b_3$}] at (0,0) {};
		\node [label={above:$a_3$}] at (0,2) {};
		\node [] at (2,0) {};
		\node [] at (2,2) {};
		\node [] at (4,0) {};
		\node [label={above:$\alpha$}] at (4,2) {};
		\node [] at (6,2) {};
		\node [] at (6,0) {};
		\node [label={below:$\beta$}] at (8,0) {};
		\node [] at (8,2) {};
		\node [] at (10,2) {};
		\node [] at (10,0) {};
	\end{scope}
		\node [label={below:$b_{3k}$}] at (2,0.05) {};
		\node [label={above:$a_{3k}$}] at (2,1.95) {};
		\node [label={below:$a_{3k+1}$}] at (4,0.1) {};
		\node [label={below:$b_{3k+1}$}] at (6,0.1) {};
		\node [label={above:$a_{3k+2}$}] at (6,1.9) {};
		\node [label={above:$b_{3k+2}$}] at (8,1.9) {};
		\node [label={above:$a_{3(k+1)}$}] at (10,1.85) {};
		\node [label={below:$b_{3(k+1)}$}] at (10,0.1) {};
\end{tikzpicture}
    \caption{The subgraph of $T_n$ whose vertical edges comprise the set $S$.}
    \label{fig:tsg-P}
\end{figure}

\begin{proof}
Label the vertices of $T_n$ as in Definition \ref{defn:trapgph}. Note that for
$i = 0, \ldots, n-1$, the vertices $a_i$ and $b_i$ each have degree 3. Without
loss of generality, consider the edge $(a_i,b_i)$; this edge belongs to exactly
two $4$-cycles, $(a_i, b_i, a_{i+1}, \alpha, a_i)$ and $(a_i, b_i, \beta, b_{i-1}, a_i)$,
where the subscripts are taken modulo $n$. Moreover, the edge $(\alpha, a_i)$
belongs to exactly two $4$-cycles, $(\alpha, a_i, b_i, a_{i+1},\alpha)$ and $(\alpha,
a_i, b_{i-1}, a_{i-1}, \alpha)$. Similarly, each edge $(\beta, b_i)$ belongs to exactly two
$4$-cycles. Since we have accounted for every edge of $T_n$, it follows that
$T_n^+$ is a $2$-regular graph. Thus, $T_n^+$ is either a cycle graph or a
disjoint union of cycle graphs.

Consider the subgraph of $T_n$ constructed by starting at the edge $(a_0,b_0)$ and consecutively finding the next (in order of increasing subscript) non-incident edge on a shared 4-cycle; a representation of this subgraph is shown in Figure \ref{fig:tsg-P}.  Let $S$ be the subset of ``vertical edges" in the figure, and denote $S^{+_{T_n}}$ by $S^+$ for simplicity.  By Lemma \ref{lem:plus-subgraph}, $S^+$ is a subgraph of $T_n^+$ and 
\begin{align*}
	V(S^+) = \{ (a_0,b_0)^+, \ &(\alpha, a_1)^+, \ (a_2,b_1)^+, \ (b_2,\beta)^+, \ (a_3,b_3)^+,\\
	& (\alpha, a_4)^+, \ (a_5,b_4)^+, \ (b_5,\beta)^+, \ (a_6,b_6)^+,\\
	& \cdots \\
	& (\alpha, a_{3k+1})^+, \ (a_{3k+2},b_{3k+1})^+, \ (b_{3k+2},\beta)^+, \ (a_{3(k+1)},b_{3(k+1)})^+, \\
	& \cdots \},
\end{align*}
where all subscripts are taken modulo $n$.

Let $j$ be the smallest positive integer such that $3j\equiv 0\pmod n$.  By construction, $S^+$ is a connected subgraph of $T_n^+$ that contains the cycle $C_m$ for $m=4j$.  So $S^+$ must equal $C_m$.  

Suppose $n\pmod 3 \not\equiv 0$.  Since $\gcd(3,n)=1$, $n$ must divide $j$ by the generalized version of Euclid's lemma.  The minimality of $j$ implies $j=n$.  Finally, since $|T_n^+| = |E(T_n)| = 4n$, we conclude that $T_n^+ \cong C_{4n}$ and $T_n^+$ is connected.  Conversely, if $T_n^+$ is connected, then $T_n^+ \cong C_{4n}$ with connected subgraph $S^+$ that contains $C_{4j}$.  So $j=n$, and $\gcd(3,n)=1$.

\end{proof}

The Clebsch graph contains the trapezohedral graph $T_3$ (the 3-cube) as an induced subgraph, and by Theorem \ref{tsgConn}, $T_3^+$ is not connected.  In fact, $\C^+$ is not connected.  Consider a four-cycle $Y$ in $\C$; using the labeling shown in Figure \ref{fig:Clebsch}, we must have $Y = (v, v^i, v^{ij}, v^j, v)$ or $Y = (v, v^i, v^{ia}, v^a, v)$, where the superscripts $i$ and $j$ indicate the the adjacent vertices differ in exactly the $i$th (respectively, $j$th) bit while the superscript $a$ indicates that the adjacent vertices differ in all four bits.  Therefore, each pair of non-incident edges of a 4-cycle must have the same bit-flip represented by the adjacency between their starting and ending vertices.  It follows that there are adjacencies only among the vertices in $\C^+$ corresponding to vertical edges in $\C$ (bit-flip in the first bit) (see Figure \ref{fig:Clebsch}), and similarly for the short diagonal edges (second bit), long diagonal edges (third bit), horizontal edges (fourth bit), and the eight extra edges not shown in the figure (all four bits).  The plus graph $\C^+$ therefore has five components, each with eight vertices and isomorphic to $K_{4,4}$.  The matrix $M$ constructed in (and preceding) the proof of Theorem \ref{thm:Clebsch-q2} assigned one of the values $\pm a, \pm b, \pm c, \pm d, \pm e$ to each of the five components of $\C^+$, but we were able to take $a=b=c=d=e$.  In general, such an assumption would be possible only if we knew that the plus graph of the given graph has a single component.  This will be our primary challenge in the next section.

\section{The Sims-Gewirtz Graph} \label{SimsGewirtz}

\subsection{Describing the Sims-Gewirtz Graph}

With 56 vertices and 280 edges, drawings of the entirety of the Sims-Gewirtz
Graph, denoted $\Gamma$ following \cite{MR1241907}, quickly become indecipherable. As such, to understand this
graph, we need a better way to visualize it. We begin with the following easy consequence of $\Gamma$ having parameters $k=10$, $\lambda = 0,$ and $\mu = 2$; similar arguments will occur in many of the proofs that follow.

\begin{lem} \label{9-4cs2}
	Any edge in $\Gamma$ is an edge in nine distinct 4-cycles that share no other edges with one another.
\end{lem}

\begin{proof}
	Let $e=(u,v)\in \Gamma$, and let $w$ be one of the nine remaining
	vertices in $N(u)$. It follows that $w\nsim v$
	since $\Gamma$ is triangle-free. Since $\mu = 2$, there exists a unique vertex $x\in N(v)$ such that $x \sim w$.  Thus, the edge $e$ belongs to the 4-cycle $(u,v,x,w,u)$.  To prove that each neighbor $w\in N(u)\setminus  \{v\}$ defines a distinct 4-cycle with $e$ as the only shared edge (that is, if $w\sim x$ and $w' \sim x'$, then $x\neq x'$), suppose $e$ belongs to $(u,v,x,w,u)$ and $(u,v,x,w',u)$.  Then $N(u)\cap N(x) = \{v,w,w'\}$, violating $\mu = 2$.  Therefore, $e$ is contained in nine distinct 4-cycles whose only shared edge is $e$.
\end{proof}

\begin{figure}[h!]
    \centering
 	\definecolor{zzttff}{rgb}{0.6,0.2,1}
\definecolor{qqzzff}{rgb}{0,0.6,1}
\definecolor{ffwwqq}{rgb}{1,0.4,0}
\definecolor{ffzzcc}{rgb}{1,0.6,0.8}
\definecolor{qqwuqq}{rgb}{0,0.39215686274509803,0}
\definecolor{yqqqqq}{rgb}{0.5019607843137255,0,0}
\definecolor{qqwwzz}{rgb}{0,0.4,0.6}
\begin{tikzpicture}[scale=1,line cap=round,line join=round,>=triangle 45,x=1cm,y=1cm]
\draw [line width=2.8pt,color=ffwwqq] (-3,3)-- (-3,-3);
\draw [line width=2.8pt,color=ffwwqq] (-1,-3)-- (-1,3);
\draw [line width=2.8pt,color=ffwwqq] (1,3)-- (1,-3);
\draw [line width=2.8pt,color=ffwwqq] (3,3)-- (3,-3);
\draw [line width=2.8pt,color=qqzzff] (-3,3)-- (3,3);
\draw [line width=2.8pt,color=qqzzff] (-3,1)-- (3,1);
\draw [line width=2.8pt,color=qqzzff] (3,-1)-- (-3,-1);
\draw [line width=2.8pt,color=qqzzff] (-3,-3)-- (3,-3);
\draw [line width=2.8pt,color=qqwuqq] (-3,3)-- (1,1);
\draw [line width=2.8pt,color=qqwuqq] (1,1)-- (3,-1);
\draw [line width=2.8pt,color=qqwuqq] (3,-1)-- (-1,-3);
\draw [line width=2.8pt,color=qqwuqq] (-1,3)-- (3,1);
\draw [line width=2.8pt,color=qqwuqq] (3,1)-- (1,-1);
\draw [line width=2.8pt,color=qqwuqq] (1,-1)-- (-3,-3);
\draw [line width=2.8pt,color=qqwuqq] (3,3)-- (-1,1);
\draw [line width=2.8pt,color=qqwuqq] (-1,1)-- (-3,-1);
\draw [line width=2.8pt,color=qqwuqq] (-3,-1)-- (1,-3);
\draw [line width=2.8pt,color=qqwuqq] (1,3)-- (-3,1);
\draw [line width=2.8pt,color=qqwuqq] (-3,1)-- (-1,-1);
\draw [line width=2.8pt,color=qqwuqq] (-1,-1)-- (3,-3);
\draw [line width=2.8pt,color=zzttff] (-3,-3)-- (-1,1);
\draw [line width=2.8pt,color=zzttff] (-1,1)-- (1,3);
\draw [line width=2.8pt,color=zzttff] (1,3)-- (3,-1);
\draw [line width=2.8pt,color=zzttff] (-3,-1)-- (-1,3);
\draw [line width=2.8pt,color=zzttff] (-1,3)-- (1,1);
\draw [line width=2.8pt,color=zzttff] (1,1)-- (3,-3);
\draw [line width=2.8pt,color=zzttff] (3,3)-- (1,-1);
\draw [line width=2.8pt,color=zzttff] (1,-1)-- (-1,-3);
\draw [line width=2.8pt,color=zzttff] (-1,-3)-- (-3,1);
\draw [line width=2.8pt,color=zzttff] (-3,3)-- (-1,-1);
\draw [line width=2.8pt,color=zzttff] (-1,-1)-- (1,-3);
\draw [line width=2.8pt,color=zzttff] (1,-3)-- (3,1);
\draw [fill=qqwwzz,draw=qqwwzz] (-1,1) circle (3.5pt);
\draw [fill=yqqqqq,draw=yqqqqq] (1,1) circle (3.5pt);
\draw [fill=qqwwzz,draw=qqwwzz] (1,-1) circle (3.5pt);
\draw [fill=yqqqqq,draw=yqqqqq] (-1,-1) circle (3.5pt);
\draw [fill=qqwwzz,draw=qqwwzz] (-3,3) circle (3.5pt);
\draw [fill=qqwuqq,draw=qqwuqq] (-1,3) circle (3.5pt);
\draw [fill=ffzzcc,draw=ffzzcc] (1,3) circle (3.5pt);
\draw [fill=yqqqqq,draw=yqqqqq] (3,3) circle (3.5pt);
\draw [fill=ffzzcc,draw=ffzzcc] (3,1) circle (3.5pt);
\draw [fill=qqwuqq,draw=qqwuqq] (3,-1) circle (3.5pt);
\draw [fill=qqwwzz,draw=qqwwzz] (3,-3) circle (3.5pt);
\draw [fill=qqwuqq,draw=qqwuqq] (1,-3) circle (3.5pt);
\draw [fill=ffzzcc,draw=ffzzcc] (-1,-3) circle (3.5pt);
\draw [fill=qqwuqq,draw=qqwuqq] (-3,1) circle (3.5pt);
\draw [fill=ffzzcc,draw=ffzzcc] (-3,-1) circle (3.5pt);
\draw [fill=yqqqqq,draw=yqqqqq] (-3,-3) circle (3.5pt);
\end{tikzpicture}
    \caption{The set $P\cup L$.}
    \label{fig:pi-design}
\end{figure}
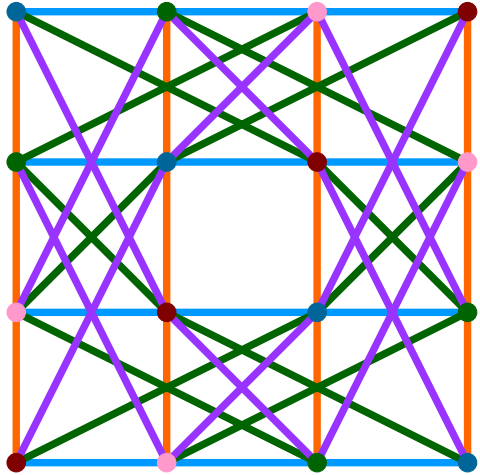

In \cite{MR1241907} Brouwer and Haemers prove that $\Gamma$ is the unique strongly regular graph with parameters $\srg(56,10,0,2)$ and give the following method for visualizing its vertices and edges.  Partition $V(\Gamma)$ into the sets ${T,P,L}$ where
$P$ and $L$ are two disjoint 16-cocliques and $T$ consists of the remaining 24 vertices.  The induced subgraph $\pi = \Gamma[P\cup L]$ is the point-block incidence graph of AG(2,4)
with one parallel class of lines removed, and the induced subgraph $\tau = \Gamma[T]$ is isomorphic to six disjoint 4-cycles \cite{MR1241907}.  It is helpful to think of the set $P$ as ``points" and the set $L$ as ``lines," where each line consists of four points, each point belongs to four lines, two non-parallel lines intersect in a unique point, and two non-parallel points are contained in a unique line.  See Figure \ref{fig:pi-design}.

Whenever two lines $l_1, l_2$ (respectively, two points $p_1, p_2$) are parallel, we denote it with $l_1 \parallel l_2$ (respectively, $p_1 \parallel p_2$).  Note that \[l_1\parallel l_2\iff N(l_1)\cap N(l_2)\subseteq T \qquad \mbox{and} \qquad p_1\parallel p_2\iff N(p_1)\cap N(p_2)\subseteq T.\]

Let $x\in P\cup L$. Since $\pi$ is 4-regular while $\Gamma$ is 10-regular, $\left|N(x)\cap T \right|=6$. Let
$t,u\in N(x)\cap T$.  Then $t,u$ cannot belong to the same 4-cycle in $\tau$: indeed, $t\nsim u$ since $\Gamma$ is triangle-free, so $t,u$ would be non-adjacent vertices in a 4-cycle in $\tau$, which would imply that
$|N(t)\cap N(u)|\geq 3$, violating $\mu=2$. So $t$ and $u$ are vertices
in disjoint 4-cycles in $\tau$. Therefore, $N(x)\cap T$ consists of one vertex
from each of the six disjoint 4-cycles in $\tau$.

Following \cite{MR1241907} (with colors replacing numbers), we label the vertices of the six 4-cycles consecutively as red, green, blue, purple; without loss of generality, the enumeration of $N(x)\cap T$ for every
$x\in P\cup L$ is shown in Figure \ref{fig:pi-labels} below.

\begin{figure}[h!]
    \centering
    \input{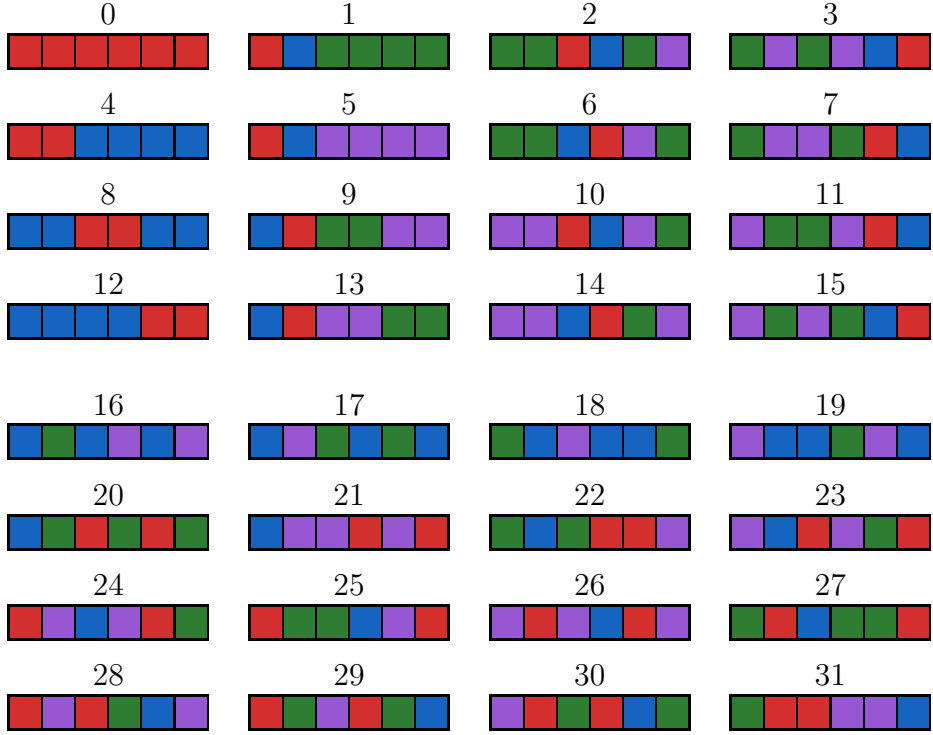}
    \caption{The 32 vertices in $\pi$, shown in terms of their neighbors in $\tau$.}
    \label{fig:pi-labels}
\end{figure}

The 16 lines in $L$ are labeled with integers $0, 1, \ldots, 15$, and the 16 points in $P$ are labeled
with integers $16, 17, \ldots, 31$ such that
\[l_1\parallel l_2\iff l_1\equiv l_2 \left({\rm mod}\: 4\right) \qquad \mbox{and} \qquad p_1\parallel p_2\iff p_1\equiv p_2 \left({\rm mod}\: 4\right).\]

Figure \ref{fig:pi-labels} describes the $32\cdot 6 = 192 = 6\cdot 4\cdot 8$ edges between the vertices of $\pi$ and the vertices of $\tau$.  After accounting for the $6\cdot 4 = 24$ edges of $\tau$, there remain the $280 - 192 - 24 = 64$ edges of the 4-regular induced subgraph graph $\pi$, which we now analyze.  We begin by considering the neighbors of a pair of non-adjacent vertices in a 4-cycle of $\tau$.

Let $t$ be an arbitrary vertex in $\tau$.  Since $\tau = 6C_4$, there exists a unique vertex $t'$ in $\tau$ such that $N(t) \cap N(t') \subseteq T$ (that is, $t$ and $t'$ belong to the same 4-cycle in $\tau$) and $t$ is not adjacent to $t'$.  In what follows, two vertices of $\tau$ denoted by the same letter with one primed will always denote such a pair.  We use this pair of vertices to define two induced subgraphs of $\Gamma$:  let
\[Z_t = \Gamma[(N(t)\cap L) \cup (N(t')\cap P) \cup \{t,t'\}]\]
    and
\[Z_{t'} = \Gamma[(N(t')\cap L) \cup (N(t)\cap P) \cup \{t',t\}].\]
We next show that the edges of these subgraphs are arranged as in Figure \ref{fig:tsg}.

\begin{lem}\label{lem:tsg}
    Given any $t\in T$, the subgraphs $Z_t$ and $Z_{t'}$ are 4-trapezohedral graphs. 
\end{lem}

\begin{proof}
    Using the description of \cite{MR1241907}, as shown in Figure \ref{fig:pi-labels}, $N(t) \cap (P\cup L)$ consists of two pairs of parallel lines and two pairs of parallel
	points.  Let $N(t) \cap L = \{l_1, l_2, m_1, m_2\}$ where $l_1\parallel l_2$ and $m_1\parallel m_2$ and $l_i \nparallel m_j$ for $i,j \in \{1,2\}$.  Since the intersection of two lines defines a unique point, for each $i,j \in \{1,2\}$, there exists a unique point $p_{ij}\in P$ such that $p_{ij} \in N(l_i)\cap N(m_j) \cap P$.  Two pairs of parallel lines must intersect in four distinct points, so the $p_{ij}$ are distinct, thereby creating the 8-cycle $(l_1, p_{11}, m_1, p_{21}, l_2, p_{22}, m_2, p_{12}, l_1)$.  Moreover, each $p_{ij}\in N(t')$ since $p_{ij} \in N(t)$ would violate $\lambda = 0$ and $p_{ij}\in N(s)\cup N(s')$ (where $\left\{s,s'\right\}=N(t)\cap N(t')$)
	would violate $\mu=2$.  So $N(t')\cap P = \{p_{11}, p_{12}, p_{21}, p_{22} \}$, and we have established that $Z_t$ is a 4-trapezohedral graph with central vertices $t$ and $t'$.  A similar argument establishes the claim for $Z_{t'}$.
\end{proof}

We will often need to refer to only the outer 8-cycle of the trapezohedral subgraphs of $\Gamma$, so we make the following definition:
\[ Z_t^{\pi} = Z_t - \{t,t'\} \quad \mbox{for } t\in T. \]

Now consider a 4-cycle $A = (t,s,t',s',t)$ in $\tau$.  We see from Figure \ref{fig:pi-labels} that the 8-cycles $Z_t^{\pi}, Z_{t'}^{\pi}, Z_s^{\pi}, Z_{s'}^{\pi}$ are mutually disjoint and contain the 32 vertices of $\pi$.  Note that any edge between $Z_t^{\pi}$ and $Z_{t'}^{\pi}$ must connect a line (respectively, point) in $N(t)$ with a point (respectively, line) in $N(t')$, but this would create a triangle in $\Gamma$.  So there are no edges in $\Gamma$ between $Z_t^{\pi}$ and $Z_{t'}^{\pi}$ or, similarly, between $Z_s^{\pi}$ and $Z_{s'}^{\pi}$.

To analyze the edges between $Z_t^{\pi}$ and $Z_s^{\pi}$ or $Z_{s'}^{\pi}$, and between $Z_{t'}^{\pi}$ and $Z_s^{\pi}$ or $Z_{s'}^{\pi}$, let $\pi_A$ be the subgraph of $\Gamma$ defined by 
\[ V(\pi_A) = V(A) \cup V(\pi) \quad\mbox{and}\quad E(\pi_A) = E(A) \cup \bigcup_{a\in A} E(Z_a). \]
Note that $\pi_A - A$ is a 2-regular spanning subgraph, containing half the edges, of $\pi$.  Starting with each of the five remaining 4-cycles in $\tau$, we similarly define $\pi_B$, $\pi_C$, $\pi_D$, $\pi_E$, and $\pi_F$.

\begin{lem} \label{lem:edge-part}
    Given a 4-cycle $A$ in $\tau$, there exists a unique 4-cycle $B$ in $\tau$ such that $\{E(\pi_A - A), E(\pi_B - B)\}$ forms a partition of $E(\pi)$.
\end{lem}

\begin{proof}
Let $A = (t,s,t',s',t)$, and as before, denote the pairs of parallel lines in
	$N(t)\cap L$ by $l_1\parallel l_2$ and $m_1\parallel m_2$.  Then there exists
	4-cycle $B = (u,v,u',v',u)$ in $\tau$ such that $N(l_1)\cap N(l_2)
	= \{t,u\}$.  In fact, it is established by Figure \ref{fig:pi-labels} that
	every pair of parallel lines  that has its first common neighbor in $A$  has its second common neighbor in $B$, and every pair of parallel points that has its first common neighbor in $A$ has its second common neighbor in $B$.  This is the pairing of 4-cycles described in \cite{MR1241907} by considering two non-parallel directions in $L$:  the parallel lines $l_1, l_2, l_3, l_4$ intersect $m_1, m_2, m_3, m_4$ in 16 points, with $l_1\parallel l_2$ intersecting $m_1\parallel m_2$ in $Z_t^{\pi}$, $l_3\parallel l_4$ intersecting $m_3\parallel m_4$ in $Z_{t'}^{\pi}$, $l_1\parallel l_2$ intersecting $m_3\parallel m_4$ in $Z_u^{\pi}$, and $l_3\parallel l_4$ intersecting $m_1\parallel m_2$ in $Z_{u'}^{\pi}$.  The other two directions correspond similarly to $Z_{s}^{\pi}$, $Z_{s'}^{\pi}$, $Z_{v}^{\pi}$, and $Z_{v'}^{\pi}$.

It is clear from the definitions of $\pi_A$ and $\pi_B$ that each vertex of $\pi$ appears twice in $V(\pi_A) \cup V(\pi_B)$.  We now claim that each edge of $\pi$ appears exactly once in $E(\pi_A) \cup E(\pi_B)$.  Indeed, suppose edge $(p,l)$ belongs to both $\pi_A$ and $\pi_B$. Without loss of generality, assume $p$ is the unique point of intersection of $l = l_i$ with $m_j$ for some $i,j\in \{1,2,3,4\}$.  Since $l_i$ and $m_j$ belong only to $Z_t^{\pi}, Z_{t'}^{\pi}, Z_u^{\pi}$, or $Z_{u'}^{\pi}$ (and not to any outer 8-cycle with subscript $s,s',v$ or $v'$), it follows that if the edge $(p,l)$ occurs twice in $E(\pi_A) \cup E(\pi_B)$, the first occurrence must be in $Z_t^{\pi}$ or $Z_{t'}^{\pi}$, and the second occurrence must be in $Z_u^{\pi}$ or $Z_{u'}^{\pi}$.  But that implies $p$ is the point of intersection of $l_i$ with $m_j$ and also the point of intersection of $l_i$ with $m_k$ for some $k \neq j$.  This contradicts the assumption that $m_j \parallel m_k$.  So $E(\pi_A-A)\cap E(\pi_B-B) = \emptyset$.  

Since $|E(\pi_A - A) \cup E(\pi_B - B)| = 64 = |E(\pi)|$, we see that every edge of the 4-regular graph $\pi$ appears exactly once in $E(\pi_A - A)\cup E(\pi_B - B)$, and the result is proved.
\end{proof}

By Lemma \ref{lem:edge-part}, we see that the six disjoint 4-cycles that comprise $\tau$ come in pairs.  For the pair of 4-cycles $A,B$ in $\tau$, define the induced subgraph $\pi_{AB}$ of $\Gamma$ by
\[ V(\pi_{AB}) = V(\pi_A) \cup V(\pi_B) = V(\pi_A) \cup V(B) = V(\pi_B) \cup V(A) \]
and
\[ E(\pi_{AB}) = E(\pi_A) \cup E(\pi_B). \]
Similarly, we define $\pi_{CD}$ and $\pi_{EF}$ using the remaining two pairs of 4-cycles in $\tau$.

\begin{figure}[h!]
    \centering
    \input{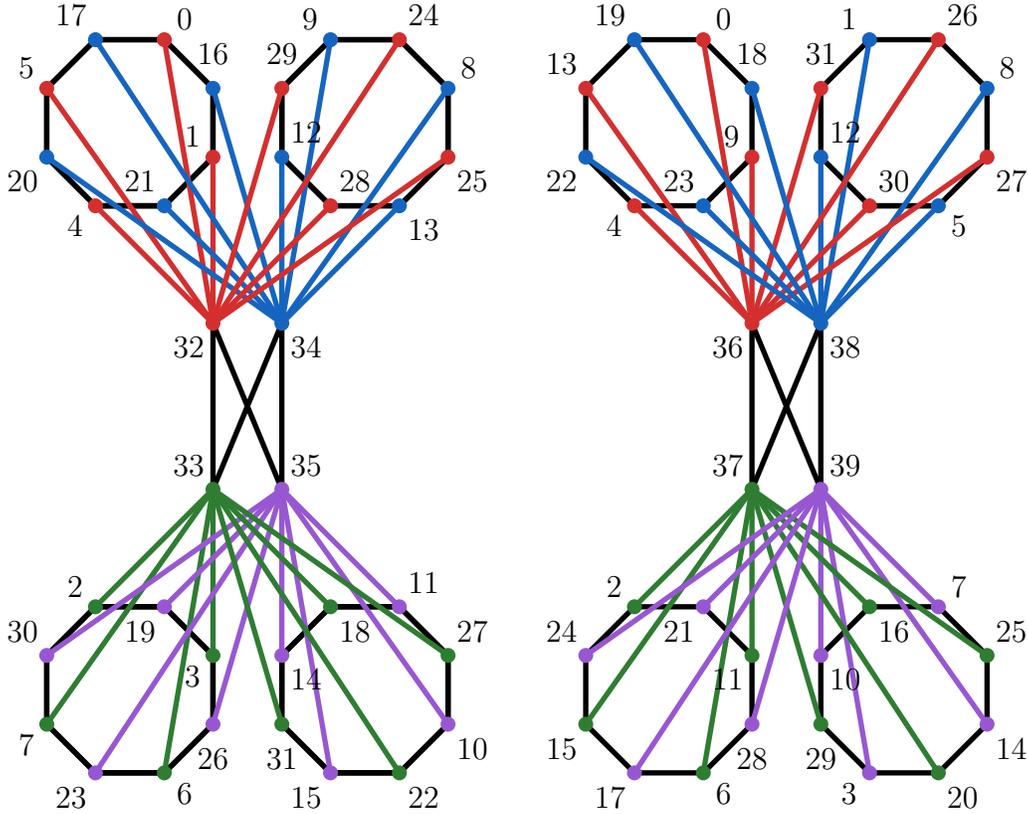}
    \caption{Subgraphs $\pi_{\mathcal{A}}$ and $\pi_{\mathcal{B}}$  comprising the induced subgraph $\pi_{\mathcal{AB}}$.}
    \label{fig:subgraph1}
\end{figure}

We use Figure \ref{fig:pi-labels} to illustrate the three partitions built by Lemma \ref{lem:edge-part}.  Starting with $\mathcal{A} = (32,33,34,35,32)$, the pairs of parallel lines are indexed by $i, i+4$ for $0\leq i \leq 3$ or $8\leq i \leq 11$, the pairs of parallel points are indexed $i, i+4$ for $16\leq i \leq 19$ or $24 \leq i \leq 27$, and $\mathcal{B} = (36,37,38,39,36)$.  Then  $\mathcal{C} = (40, 41, 42, 43, 40)$ and $\mathcal{D} = (44,45,46,47,44)$ are the pair of 4-cycles corresponding to the pairs of parallel lines $i,i+8$ for $0\leq i \leq 7$ and the pairs of parallel points $i,i+8$ for $16\leq i \leq 23$.  Finally, $\mathcal{E} = (48,49,50,51,48)$ and $\mathcal{F} = (52,53,54,55,52)$ correspond to the pairs of parallel lines $i,i+12$ for $0\leq i \leq 3$ or $8\leq i \leq 11$ (the latter subtracted by 16) and the pairs of parallel points $i, i+12$ for $16\leq i \leq 19$ or $24\leq i \leq 27$ (the latter subtracted by 16).  

\begin{figure}[h!]
    \centering
    \input{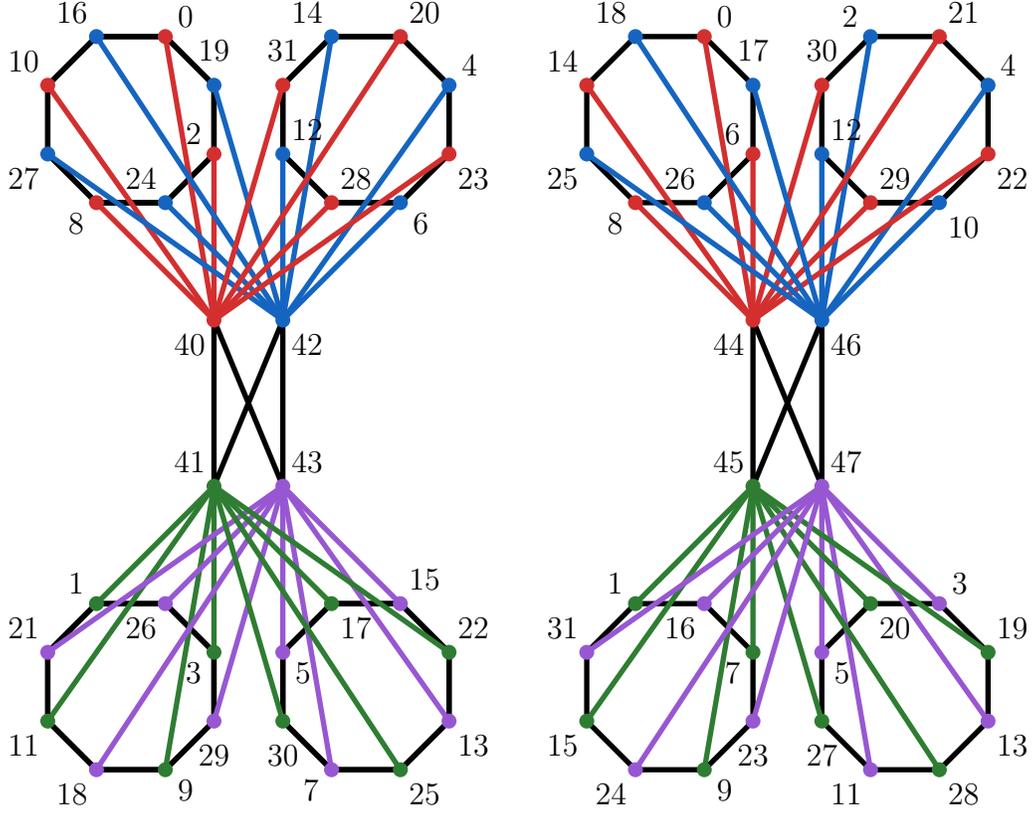}
    \caption{Subgraphs $\pi_{\mathcal{C}}$ and $\pi_{\mathcal{D}}$  comprising the induced subgraph $\pi_{\mathcal{CD}}$.}
    \label{fig:subgraph2}
\end{figure}

\begin{figure}[h!]
    \centering
    \input{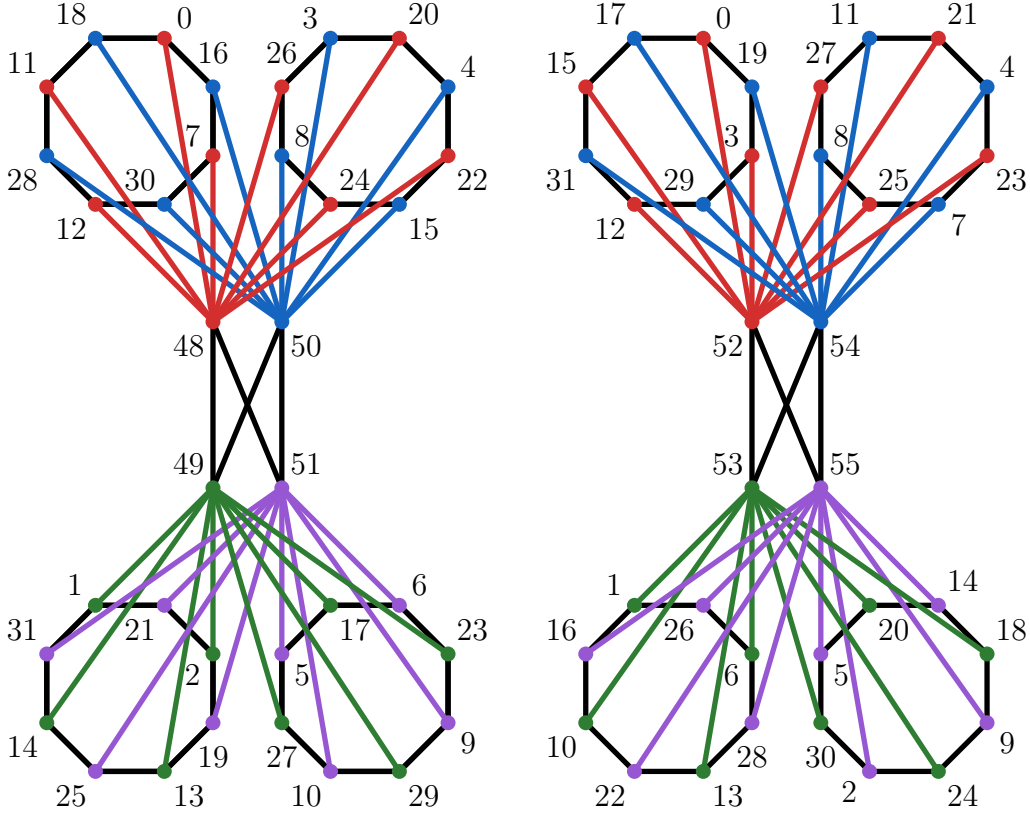}
    \caption{Subgraphs $\pi_{\mathcal{E}}$ and $\pi_{\mathcal{F}}$ comprising the induced subgraph $\pi_{\mathcal{EF}}$.}
    \label{fig:subgraph3}
\end{figure}

The three induced subgraphs $\pi_{\mathcal{AB}}$, $\pi_{\mathcal{CD}}$, and $\pi_{\mathcal{EF}}$ together show $\Gamma$ in its entirety:  each vertex of $\pi$ appears six times, each edge of $\pi$ appears three times (once per 4-cycle pairing), each vertex of $\tau$ appears once, and each edge in $E(\Gamma)\setminus E(\pi)$ appears once. See Figures \ref{fig:subgraph1}, \ref{fig:subgraph2}, and \ref{fig:subgraph3}.

\subsection{The Plus Graph of the Sims-Gewirtz Graph}

Now that we have a full view of the Sims-Gewirtz graph $\Gamma$, we will prove that its plus graph $\Gamma^+$ is connected.  We first show that the edges in the outer 8-cycle of any single trapezohedral subgraph in one of the induced subgraphs defined by a pairing of two 4-cycles in $\tau$ are fully separated among the eight trapezohedral subgraphs defined by another pairing.

\begin{lem}\label{C-covers2}
	Let $A,B,C,D$ be four 4-cycles in $\tau$ such that $A$ is paired with $B$ and $C$
	is paired with $D$. If two distinct edges $(p,l), (q,m) \in
	E(Z^\pi_a)$ for some $a\in V(A)\cup V(B)$, then there exists no $c\in V(C)\cup V(D)$ such that $(p,l), (q,m) \in E(Z_c^{\pi})$.  
\end{lem}

\begin{proof}
    Suppose $(p,l), (q,m) \in E(Z_a^{\pi})$ for some $a\in V(A)\cup V(B)$.  Without loss of generality, we have three cases to consider:  $l = m$, and if $l \neq m$, then $l\parallel m$ or $l\nparallel m$.  If $l = m$, then $N(p) \cap N(q) = \{l, a'\}$, so that $\mu=2$ prevents $p$ and $q$ also sharing $c'$ as a common neighbor.  Suppose $l\neq m$.  If $l \parallel m$, then $N(l) \cap N(m) = \{a,b\}$, where $b\in B$, and again $\mu = 2$ prevents $l$ and $m$ also sharing $c$ as a common neighbor.  Finally, if $l \nparallel m$, then they share two common neighbors in $Z_a$ and again cannot both be adjacent to $c$.
\end{proof}

\begin{cor} \label{tsg-ring-pm}
	Let $A,B,C$ be 4-cycles in $\tau$ such that $A$ and $B$ are paired with each other and not paired with $C$, and let $c\in C$. The sets
	$E(Z_c^\pi)\cap E(\pi_A)$ and $E(Z_c^\pi)\cap E(\pi_B)$ are two perfect matchings of $Z_c^\pi$.
\end{cor}

\begin{proof}
    Let $e$ and $f$ be two incident edges in $Z_c^{\pi}$.  Without loss of generality, assume $e = (p,l)$ and $f = (p,m)$ for some $p\in P$ and $l,m \in L$.  By Lemmas \ref{lem:edge-part} and \ref{C-covers2}, $e$ and $f$ belong to $E(\pi_{AB})$ but do not belong to the outer 8-cycle of any one trapezohedral subgraph in $\pi_{AB}$.  Suppose, without loss of generality, that $e \in E(Z_a^{\pi})$ for some $a\in A$.  Since $f\notin E(Z_a^{\pi})$, $f\notin E(\pi_A)$ as $p$ appears exactly once in $\pi_A$.  By Lemma \ref{lem:edge-part}, $f\in E(\pi_B)$.  The result follows by noting that the eight edges in $Z_c^{\pi}$ therefore alternate between belonging to $\pi_A$ and $\pi_B$.
\end{proof}

Recall from Theorem \ref{tsgConn} that the plus graph of any individual 4-trapezohedral subgraph is connected.  This is the key fact needed to prove the main result of this subsection.

\begin{thm} \label{thm:connected}
	$\Gamma^+$ is connected.
\end{thm}

\begin{proof}
We decompose the edge set of $\Gamma$ as $E(\Gamma) = E(\pi) \cup E(\tau) \cup S$ where $S = E(\Gamma) \setminus (E(\pi)\cup E(\tau))$.  Since all plus graphs considered will be with respect to $\Gamma$, we will omit the $\Gamma$ subscript on the superscript $^+$ for any subgraph or set of edges of $\Gamma$.  First we show that $\pi^+$ is connected.

    Let $e^+$ and $f^+$ be two distinct vertices of $\pi^+$.  Then $e,f \in E(\pi_{AB}) \cap E(\pi)$ for some paired 4-cycles $A,B$ in $\tau$ by Lemma \ref{lem:edge-part}.  If there exists $a\in V(A)\cup V(B)$ such that $e,f \in Z_a^{\pi}$, then there exists a path in $Z_a^+$ from $e^+$ to $f^+$ by Theorem \ref{tsgConn}.  Otherwise, suppose $e\in E(Z_a^{\pi})$ and $f \in E(Z_b^{\pi})$ for some $a\neq b \in V(A) \cup V(B)$.  Considering these edges in $\pi_{CD}$ for paired 4-cycles $C,D$ in $\tau$, there exists $c,d \in V(C) \cup V(D)$ such that $e \in E(Z_c^{\pi})$ and $f \in E(Z_d^{\pi})$.  If $c=d$, then there exists a path in $Z_c^+$ from $e^+$ to $f^+$ by Theorem \ref{tsgConn}.  Suppose $c\neq d$.  By Lemmas \ref{lem:edge-part} and \ref{C-covers2}, the edges of $Z_a^{\pi}$ are distributed so that each edge appears in exactly one trapezohedral subgraph of $\pi_{CD}$.  So there exists some edge $h \in E(Z_a^{\pi})$ with $h \in E(Z_d^{\pi})$.  Then by applying Theorem \ref{tsgConn} twice, we can find a path in $Z_a^+$ from $e^+$ to $h^+$ and a path in $Z_d^+$ from $h^+$ to $f^+$.

    The connectedness of $\Gamma^+$ will now follow if, starting at any vertex of $(E(\tau)\cup S)^+$, we can establish the existence of a path ending at a vertex of $\pi^+$.  Given any edge $(w,x) \in S$ with $w \in V(\tau)$ and $x\in V(\pi)$, we have $(w,x)^+ \in Z_w^+$ or $(w,x)^+ \in Z_{w'}^+$. Theorem \ref{tsgConn} immediately provides a path from $(w,x)^+$ to a vertex of $\pi^+$.

    Finally, by Lemma \ref{9-4cs2}, any edge $(t,s) \in E(\tau)$ is in nine 4-cycles.
	However, only one of these 4-cycles is completely in $\tau = 6C_4$, so eight must each be of the form $(t,s,p,l,t)$ or $(t,s,l,p,t)$ for some $p\in P$ and $l\in L$.  Then $(p,l)^+ \in V(\pi^+)$ is such that  $(t,s)^+ \sim (p,l)^+$, completing the proof.
 \end{proof}

\subsection{Proof of the Sims-Gewirtz Graph Not Having Two Eigenvalues.}

In this subsection we prove the impossibility of the existence of a Gram matrix $M \in \sm(\Gamma)$ such that $M^2 = M$, thereby establishing that $q(\Gamma) \neq 2$ and hence $q(\Gamma) = 3$. 

\begin{lem} \label{lem:abs-val}
    Suppose $q(\Gamma)=2$, and let $M \in
	\sm(\Gamma)$ be such that $M = [ \langle f_i, f_j \rangle ]$ and $M^2 = M$.  Then all nonzero off-diagonal entries of $M$ have the same absolute value.
\end{lem}

\begin{proof}
    By Lemma \ref{signlemma}, $M_{ij} = \pm M_{kl}$ whenever $(v_i, v_j, v_k, v_l, v_i)$ is a 4-cycle in $\Gamma$.  By Theorem \ref{thm:connected}, $\Gamma^+$ is connected, and the result follows.
\end{proof} 

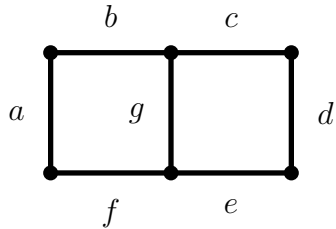
\begin{figure}[h!]
		\centering
		\begin{tikzpicture}[scale=0.8, line cap=round,line join=round,>=triangle 45,x=1cm,y=1cm,V/.style={circle, inner sep = 2pt,fill=black}]
\draw [line width=2pt] (-2,0)-- (-2,2) node [midway, label={left:$a$}] {};
\draw [line width=2pt] (-2,2)-- (0,2) node [midway, label={above:$b$}] {};
\draw [line width=2pt] (0,2)-- (2,2) node [midway, label={above:$c$}] {};
\draw [line width=2pt] (2,2)-- (2,0) node [midway, label={right:$d$}] {};
\draw [line width=2pt] (2,0)-- (0,0) node [midway, label={below:$e$}] {};
\draw [line width=2pt] (0,0)-- (-2,0) node [midway, label={below:$f$}] {};
\draw [line width=2pt] (0,0)-- (0,2) node [midway, label={left:$g$}] {};

\begin{scope}[nodes=V]
\node at (-2,0) {};
\node at (-2,2) {};
\node at (2,2) {};
\node at (2,0) {};
\node at (0,2) {};
\node at (0,0) {};
\end{scope}
\end{tikzpicture}
		\caption{Crossbar 6-cycle.}
		\label{fig:crossbar}
\end{figure}

Due to Lemma \ref{lem:abs-val}, we will care only about the sign of each nonzero off-diagonal entry $M_{ij}$, that is, the sign of each edge.  Through an abuse of notation, we will use single letters to represent both the edge $(v_i, v_j)$ and the associated number $\mathrm{sgn}(M_{ij}) = \pm 1$.  Often it will be convenient to write cycles in terms of edges instead of vertices, and we will enclose them in brackets instead of parentheses. For example, we define a {\em crossbar 6-cycle} as the sequence of edges $[a,b,c,d,e,f]$ (where $a$ is incident with $f$) such that there exists an additional edge $g$ that is incident with $b,c,e,f$, as in Figure \ref{fig:crossbar}.

\begin{defn}
	Given a cycle $C$ in $\Gamma$, we call $C$ an {\em odd} cycle if
	\[\prod_{e\in E(C)} e=-1\]
	and {\em even} if
	\[\prod_{e\in E(C)} e=1.\]
\end{defn}

\begin{lem}\label{odd4c}
	If $q(\Gamma)=2$, then all 4-cycles in $\Gamma$ are odd.
\end{lem}

\begin{proof}
    Let $M \in \sm(\Gamma)$ be such that $M = [ \langle f_i, f_j \rangle ]$ and $M^2 = M$.  By Lemma \ref{signlemma}, $M_{ij} = \pm M_{kl}$ and $M_{jk} = \mp M_{il}$, with exactly one choice positive and one negative, whenever $(v_i, v_j, v_k, v_l, v_i)$ is a 4-cycle in $\Gamma$.
\end{proof}

\begin{lem}\label{crossbar}
	If $q(\Gamma) = 2$, then all crossbar 6-cycles in $\Gamma$ are even.
\end{lem}

\begin{proof}
	Let $C_6 = [a,b,c,d,e,f]$ in Figure \ref{fig:crossbar}. Since
	$abgf=gcde=-1$
	by Lemma \ref{odd4c},
	\begin{align*}
		abcdef = (abf)(cde) = (-g)(-g) = 1.
	\end{align*}
\end{proof}

\begin{prop}\label{even-8-cyc}
	Suppose $q(\Gamma) = 2$, and let $C_8$ be an 8-cycle and $v$ be a vertex in $\Gamma$.  If $\left|N(v)\cap
	V(C_8)\right|=4$, then $C_8$ is even.
\end{prop}

\begin{proof}
	Let $C_8=[a,b,c,d,e,f,g,h]$ and let the edges incident to $v$ be $i,j,k,l$.  Without loss of generality, suppose edge $i$ is incident with the pair of edges $\{a,b\}$.  Since $\lambda=0$, it follows that each edge $j,k,l$ is incident with exactly one pair of incident edges $\{c,d\},\{e,f\},\{g,h\}$, and we may assume the correspondence is in the order listed here.  Then $[h,a,b,c,j,l]$ and $[d,e,f,g,l,j]$ are crossbar 6-cycles, so by Lemma \ref{crossbar},
	\[abcdefgh=(habc)(defg)=(jl)(jl)=1.\]
\end{proof}

Since any trapezohedral subgraph $T_4$ in $\Gamma$ satisfies the hypotheses of the previous lemma, we have the following immediate corollary.

\begin{cor} \label{even-TSG}
	If $q(\Gamma) = 2$, then $Z_t^\pi$ is an even 8-cycle for any $t\in V(\tau)$.
\end{cor}

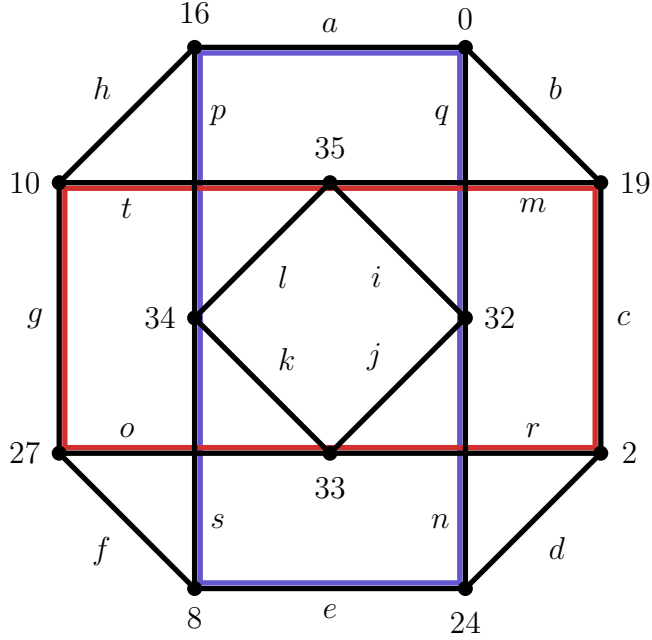
\begin{figure}[h]
	\centering
\definecolor{R}{rgb}{0.8274509803921568,0.1843137254901961,0.1843137254901961}
\definecolor{P}{rgb}{0.396078431372549,0.3411764705882353,0.8235294117647058}
\begin{tikzpicture}[scale=.6,line cap=round,line join=miter,>=triangle 45,x=1cm,y=1cm,
	every node/.style={inner sep = 2pt},
	V/.style={circle, fill=black}]

	%
	%
	%

	\draw [line width=2pt,color=P] (-2.875,5.875)-- (2.857,5.875);
	\draw [line width=2pt,color=P] (-2.875,5.875)-- (-2.875,-5.875);
	\draw [line width=2pt,color=P] (-2.875,-5.875)-- (2.875,-5.875);
	\draw [line width=2pt,color=P] (2.875,-5.875)-- (2.875,5.875);

	\draw [line width=2pt,color=R] (5.875,-2.875)-- (5.875,2.857);
	\draw [line width=2pt,color=R] (5.875,-2.875)-- (-5.875,-2.875);
	\draw [line width=2pt,color=R] (-5.875,-2.875)-- (-5.875,2.875);
	\draw [line width=2pt,color=R] (-5.875,2.875)-- (5.875,2.875);

	\draw [line width=2pt] (3,6)-- (-3,6) node [midway, label={above:$a$}] {};
	\draw [line width=2pt] (6,3)-- (3,6) node [midway, label={above right:$b$}] {};
	\draw [line width=2pt] (6,-3)-- (6,3) node [midway, label={right:$c$}] {};
	\draw [line width=2pt] (3,-6)-- (6,-3) node [midway, label={below right:$d$}] {};
	\draw [line width=2pt] (-3,-6)-- (3,-6) node [midway, label={below:$e$}] {};
	\draw [line width=2pt] (-6,-3)-- (-3,-6) node [midway, label={below left:$f$}] {};
	\draw [line width=2pt] (-6,3)-- (-6,-3) node [midway, label={left:$g$}] {};
	\draw [line width=2pt] (-3,6)-- (-6,3) node [midway, label={above left:$h$}] {};

	\draw [line width=2pt] (3,0)-- (0,3) node [midway, label={below left:$i$}] {};
	\draw [line width=2pt] (0,-3)-- (3,0) node [midway, label={above left:$j$}] {};
	\draw [line width=2pt] (-3,0)-- (0,-3) node [midway, label={above right:$k$}]{};
	\draw [line width=2pt] (0,3)-- (-3,0) node [midway, label={below right:$l$}] {};

	\draw [line width=2pt] (6,3)-- (0,3) node [pos=.25, label={below:$m$}] {};
	\draw [line width=2pt] (3,-6)-- (3,0) node [pos=.25, label={left:$n$}] {};
	\draw [line width=2pt] (-6,-3)-- (0,-3) node [pos=.25, label={above:$o$}] {};
	\draw [line width=2pt] (-3,6)-- (-3,0) node [pos=.25, label={right:$p$}] {};

	\draw [line width=2pt] (3,0)-- (3,6) node [pos=.75, label={left:$q$}] {};
	\draw [line width=2pt] (0,-3)-- (6,-3) node [pos=.75, label={above:$r$}] {};
	\draw [line width=2pt] (-3,0)-- (-3,-6) node [pos=.75, label={right:$s$}] {};
	\draw [line width=2pt] (0,3)-- (-6,3) node [pos=.75, label={below:$t$}] {};

	\begin{scope}[nodes=V]
		\node [label={above:16}] at (-3,6) {};
		\node [label={above:0}] at (3,6) {};
		\node [label={right:19}] at (6,3) {};
		\node [label={right:2}] at (6,-3) {};
		\node [label={below:24}] at (3,-6) {};
		\node [label={below:8}] at (-3,-6) {};
		\node [label={left:27}] at (-6,-3) {};
		\node [label={left:10}] at (-6,3) {};
		\node [label={above:35}] at (0,3) {};
		\node [label={right:32}] at (3,0) {};
		\node [label={below:33}] at (0,-3) {};
		\node [label={left:34}] at (-3,0) {};
	\end{scope}
\end{tikzpicture}
	\caption{Induced subgraph on $V(\mathcal{A})\cup V(Z_{40}^\pi)$ with edges
	labeled.}
	\label{fig:cross-induced-graph}
\end{figure}

\begin{lem} \label{even-purple}
    Suppose $q(\Gamma) = 2$.  For the 4-cycles $\mathcal{A} = (32,33,34,35,32)$ and $\mathcal{C} = (40,41,42,43,40)$ in $\tau$, the induced subgraph $\Gamma[V(Z_{40}^\pi)\cup V(\mathcal{A})]$ appears as in Figure \ref{fig:cross-induced-graph}.  Given the edge labels as seen in the figure, the 6-cycle $[a,q,n,e,s,p]$ (indicated in blue) or the 6-cycle $[g,t,m,c,r,o]$ (indicated in red) must be even.
\end{lem}

\begin{proof}
The first claim can be established by observation of Figures \ref{fig:subgraph1} and \ref{fig:subgraph2}.  (In fact, given any two non-paired 4-cycles $A,C$ in $\tau$, the subgraph induced by the vertices of $Z_u^\pi$, where $u\in C$, and $A$ must be isomorphic to the graph shown in Figure \ref{fig:cross-induced-graph}, as a consequence of the $\lambda = 0$ and $\mu = 2$ restrictions.)
	Suppose both $[a,q,n,e,s,p]$ and $[g,t,m,c,r,o]$ are
	odd 6-cycles. Then it follows that $ae=-qnsp$ and $gc=-tmro$. Additionally,
	by Lemma \ref{odd4c},
	$b=-qim, d=-njr, f=-sko, h=-plt$. Thus,
	\begin{align*}
		abcdefgh &= (ae)(gc)(b)(d)(f)(h) \\
		&= (-qnsp)(-tmro)(-qim)(-njr)(-sko)(-plt) \\
		&= q^2n^2s^2p^2t^2m^2r^2o^2ijkl \\
		&= ijkl \\
		&= -1
	\end{align*}
	However, this contradicts the fact that $abcdefgh=1$ by Corollary
	\ref{even-TSG}.
\end{proof}

We are ready to prove our main result.

\begin{figure}[h]
	\centering
\definecolor{R}{rgb}{0.8274509803921568,0.1843137254901961,0.1843137254901961}
\definecolor{G}{rgb}{0.1803921568627451,0.49019607843137253,0.19607843137254902}
\definecolor{B}{rgb}{0.08235294117647059,0.396078431372549,0.7529411764705882}
\definecolor{P}{rgb}{0.396078431372549,0.3411764705882353,0.8235294117647058}
\definecolor{O}{rgb}{.753,.082,.396}
\begin{tikzpicture}[scale=1,line join=round,>=triangle 45,x=1cm,y=1cm,
	every node/.style={inner sep = 2pt},
	V/.style={circle}]

	\draw [line width=2pt] (4,6)-- (-4,6) node [midway, label={above right:$a$}] {};
	\draw [line width=3pt,color=R] (6,4)-- (4,6) node [midway, label={above right:$b$}] {};
	\draw [line width=3pt,color=R] (6,-4)-- (6,4) node [midway, label={right:$c$}] {};
	\draw [line width=3pt,dash pattern={on 10pt off 10pt},color=B] (6,-4)-- (6,4) node [midway, label={right:$c$}] {};
	\draw [line width=3pt,color=B] (4,-6)-- (6,-4) node [midway, label={below right:$d$}] {};
	\draw [line width=2pt] (-4,-6)-- (4,-6) node [midway, label={below right:$e$}] {};
	\draw [line width=2pt] (-6,-4)-- (-4,-6) node [midway, label={below left:$f$}] {};
	\draw [line width=2pt] (-6,4)-- (-6,-4) node [midway, label={left:$g$}] {};
	\draw [line width=2pt] (-4,6)-- (-6,4) node [midway, label={above left:$h$}] {};

	\draw [line width=2pt] (4,0)-- (0,4) node [midway, label={above right:$i$}] {};
	\draw [line width=2pt] (0,-4)-- (4,0) node [midway, label={below right:$j$}] {};
	\draw [line width=2pt] (-4,0)-- (0,-4) node [midway, label={below left:$k$}]{};
	\draw [line width=2pt] (0,4)-- (-4,0) node [midway, label={above left:$l$}] {};

	\draw [line width=3pt,color=B] (6,4)-- (0,4) node [pos=.2, label={below:$m$}] {};
	\draw [line width=3pt,color=B] (4,-6)-- (4,0) node [pos=.2, label={left:$n$}] {};
	\draw [line width=2pt] (-6,-4)-- (0,-4) node [pos=.2, label={above:$o$}] {};
	\draw [line width=2pt] (-4,6)-- (-4,0) node [pos=.2, label={right:$p$}] {};

	\draw [line width=3pt,color=R] (4,0)-- (4,6) node [pos=.8, label={left:$q$}] {};
	\draw [line width=3pt,color=R] (0,-4)-- (6,-4) node [pos=.8, label={above:$r$}] {};
	\draw [line width=2pt] (-4,0)-- (-4,-6) node [pos=.8, label={right:$s$}] {};
	\draw [line width=2pt] (0,4)-- (-6,4) node [pos=.8, label={below:$t$}] {};

	\draw [line width=3pt,color=B] ( .5, 1.5)-- (-.5, 1.5) node [midway,label={below:$a'$}] {};
	\draw [line width=2pt] ( 1.5, .5)-- ( .5, 1.5) node [midway,label={below left:$b'$}] {};
	\draw [line width=2pt] ( 1.5,-.5)-- ( 1.5, .5) node [midway,label={left:$c'$}] {};
	\draw [line width=2pt] ( .5,-1.5)-- ( 1.5,-.5) node [midway,label={above left:$d'$}] {};
	\draw [line width=3pt,color=R] (-.5,-1.5)-- ( .5,-1.5) node [midway,label={above:$e'$}] {};
	\draw [line width=3pt,color=R] (-1.5,-.5)-- (-.5,-1.5) node [midway,label={above right:$f'$}] {};
	\draw [line width=2pt] (-1.5, .5)-- (-1.5,-.5) node [midway,label={right:$g'$}] {};
	\draw [line width=3pt,color=B] (-.5, 1.5)-- (-1.5, .5) node [midway,label={below right:$h'$}] {};
	\draw [line width=2pt] ( 1.5, .5)-- ( 0.0, 4.0) node [midway,label={right:$m'$}] {};
	\draw [line width=3pt,color=R] ( .5,-1.5)-- ( 4.0, 0.0) node [midway,label={below:$n'$}] {};
	\draw [line width=3pt,color=R] (-1.5,-.5)-- ( 0.0,-4.0) node [midway,label={left:$o'$}] {};
	\draw [line width=2pt] (-.5, 1.5)-- (-4.0, 0.0) node [midway,label={above:$p'$}] {};
	\draw [line width=3pt,color=B] ( 4.0, 0.0)-- ( .5, 1.5) node [midway,label={above:$q'$}] {};
	\draw [line width=2pt] ( 0.0,-4.0)-- ( 1.5,-.5) node [midway,label={right:$r'$}] {};
	\draw [line width=2pt] (-.5,-1.5)-- (-4.0, 0.0) node [midway,label={below:$s'$}] {};
	\draw [line width=3pt,color=B] ( 0.0, 4.0)-- (-1.5, .5) node [midway,label={left:$t'$}] {};

	\begin{scope}[nodes=V]
		\node [fill=black,label={above:16}] at (-4,6) {};
		\node [fill=black,label={above:0}] at (4,6) {};
		\node [fill=black,label={right:19}] at (6,4) {};
		\node [fill=black,label={right:2}] at (6,-4) {};
		\node [fill=black,label={below:24}] at (4,-6) {};
		\node [fill=black,label={below:8}] at (-4,-6) {};
		\node [fill=black,label={left:27}] at (-6,-4) {};
		\node [fill=black,label={left:10}] at (-6,4) {};
		\node [fill=black,label={above:35}] at (0,4) {};
		\node [fill=black,label={right:32}] at (4,0) {};
		\node [fill=black,label={below:33}] at (0,-4) {};
		\node [fill=black,label={left:34}] at (-4,0) {};
		\node [fill=black,label={right:23}] at (1.5,.5) {};
		\node [fill=black,label={right:6}] at (1.5,-.5) {};
		\node [fill=black,label={left:14}] at (-1.5,.5) {};
		\node [fill=black,label={left:31}] at (-1.5,-.5) {};
		\node [fill=black,label={above:4}] at (.5,1.5) {};
		\node [fill=black,label={above:20}] at (-.5,1.5) {};
		\node [fill=black,label={below:12}] at (-.5,-1.5) {};
		\node [fill=black,label={below:28}] at (.5,-1.5) {};
	\end{scope}
\end{tikzpicture}
	\caption{Induced subgraph on $V(\mathcal{A})\cup V(Z_{40}^\pi)\cup V(Z_{42}^\pi)$ with edges labeled.}
	\label{fig:final-graph}
\end{figure}

\begin{thm}
  For the Sims-Gewirtz graph $\Gamma$, $q(\Gamma) \neq 2$ and hence $q(\Gamma) = 3$.
\end{thm}

\begin{proof}
Suppose to the contrary that $q(\Gamma) = 2$.  Using the 4-cycles $\mathcal{A} = (32,33,34,35,32)$ and $\mathcal{C} = (40,41,42,43,40)$, we perform the construction in Lemma $\ref{even-purple}$ for $Z_{40}^{\pi}$ (with edge labels $a$---$h$ and $m$---$t$) and then repeat for $Z_{42}^{\pi}$ (with edge lables $a'$---$h'$ and $m'$---$t'$).  Since there are no edges between the vertices of $Z_{40}^{\pi}$ and those of $Z_{42}^{\pi}$, the induced subgraph $\Gamma[V(\mathcal{A})\cup V(Z_{40}^\pi)\cup V(Z_{42}^\pi)]$ appears as in Figure \ref{fig:final-graph}.

	By Lemma \ref{even-purple}, the 6-cycles $[p',a',q',n',e',s']$ and
	$[m',c',r',o',g',t']$ cannot both be odd. Suppose that the first is even.  Then consider the 8-cycle $C_8 = [c,b,q,n',e',f',o',r]$ (indicated in red in Figure \ref{fig:final-graph}).
	This 8-cycle is even by Proposition \ref{even-8-cyc} since $N(40) \cap V(C_8) = \{0,2,28,31\}$.
    Repeated application of Lemma \ref{odd4c} implies
    \begin{align*}
		1&=c(bq)n'e'(f'o')r=c(-mi)n'e'(-s'k)r=(ik)cmn'e's'r \\
		&=(-jl)cmn'e's'r \\
        &=-jlcm(n'e's')r = -jlcm(p'a'q')r  \\
        &= -cm(lp')a'q'(jr) = -cm(-t'h')a'q'(-nd) \\
       &= -cdnq'a'h't'm.
    \end{align*}    
However, the 8-cycle $D_8 = [c,d,n,q',a',h',t',m]$ (indicated in blue in Figure \ref{fig:final-graph}) is even by Proposition \ref{even-8-cyc} since $N(42) \cap V(D_8) = \{4,14,19,24\}$, yielding a contradiction.

On the other hand, if the 6-cycle $[m',c',r',o',g',t']$ is even, then we begin at the 8-cycle $[e,d,r,o',g',h',p',s]$ and end at the 8-cycle $[e,f,o,r',c',b',q',n]$; arguing as above, we obtain the contradiction 
\[ 1 = edro'g'h'p's = -efor'c'b'q'n = -1.\]
\end{proof}

Having shown that, among the first six triangle-free strongly regular graphs, the Clebsch graph is unique in describing a matrix with exactly two distinct eigenvalues, we leave open this question for the seventh, Higman-Sims graph, or any triangle-free SRG yet to be discovered.

\section*{Acknowledgements}

For preliminary exploratory analysis, this work used Anvil at Purdue University through allocation \#MTH240035 from the Advanced Cyberinfrastructure Coordination Ecosystem: Services \& Support (ACCESS) program \cite{access}, which is supported by U.S. National Science Foundation grants \#2138259, \#2138286, \#2138307, \#2137603, and \#2138296.  The authors also gratefully acknowledge support in part by the NSF through grant DMS-2331072.  The results in Subsection \ref{q(Clebsch)=2} were obtained jointly with Dayton Singer and supported in part by a Fort Lewis College Betty Haskell Mathematics Scholars Grant.


\nocite{*}
\bibliographystyle{plainurl}
\bibliography{main.bib}

\end{document}